%% file: sbm_DG_mg.tex
\documentclass{elsarticle}

\usepackage{amsmath,amssymb,amsfonts,graphicx,hyperref}
\usepackage[margin=1in]{geometry}

\input{preamble.tex}

\makeatletter
\def\ps@pprintTitle{%
 \let\@oddhead\@empty
 \let\@evenhead\@empty
 \def\@oddfoot{\reset@font\hfil\thepage\hfil}%
 \let\@evenfoot\@oddfoot}
\makeatother

\title{A Geometric Multigrid Preconditioner for Discontinuous Galerkin Shifted Boundary Method}
\author{Micha\l{} Wichrowski\footnote{Heidelberg University, Germany, mt.wichrowsk@uw.edu.pl}}

\begin{document}

\begin{abstract}
    This paper introduces a geometric multigrid preconditioner for the Shifted Boundary Method (SBM) designed to solve
    PDEs on complex geometries. While the SBM simplifies mesh generation by using a non-conforming background grid, it often results in non-symmetric linear systems, a complete analysis of which is missing from the literature, especially with regard to the algebraic solution and preconditioning of high-order versions of SBM.
    Standard multigrid methods with pointwise smoothers prove ineffective for such systems due to the localized
    perturbations introduced by the shifted boundary conditions. To address this challenge, we introduce a
    Discontinuous
    Galerkin (DG) formulation for SBM that enables the design of a cell-wise multiplicative smoother within an
    $hp$-multigrid framework. The element-local nature of DG methods naturally facilitates cell-wise correction, which
    can effectively handle the local complexities arising from the boundary treatment. Numerical results for the
    Poisson equation demonstrate favorable performance with mesh refinement for linear ($p=1$) and quadratic ($p=2$)
    elements in both 2D and 3D, with iteration counts showing mild growth.
    However, challenges emerge for cubic ($p=3$) elements, particularly in 3D, where the current smoother shows
    reduced effectiveness.
\end{abstract}
\maketitle
\section{Introduction}

Solving partial differential equations on domains with complex or evolving geometries presents a significant challenge.
Traditional numerical methods, such as the finite element method (FEM), often require body-fitted meshes, which can be
difficult to generate for intricate shapes. Unfitted finite element methods offer an attractive alternative by
employing a fixed background mesh that does not conform to the physical domain boundary. Among these, the Shifted
Boundary Method (SBM)~\cite{main2018shifted} stands out by defining the problem on a surrogate domain composed of cells
from the background mesh and extrapolating boundary conditions from the true boundary to the boundary of these selected
cells. This approach avoids complex cell-cutting procedures and specialized quadrature rules common in methods like
CutFEM~\cite{burman2015cutfem}. However, the resulting linear system still has a condition number scaling like
$\mathcal{O}(h^{-2})$\cite{atallah2021shifted}, making it challenging for iterative solvers. Our results for linear
elements ($p=1$) confirm the findings in~\cite{atallah2021shifted, collins2023penalty} that SBM exhibits conditioning
similar to body-fitted methods. The performance of algebraic solvers associated with the SBM has not been thoroughly
studied and the preconditioning and efficient numerical solution of the algebraic systems emanating from SBM
formulations has remained somewhat unexplored, particularly in the context of high-order methods. This work is
attempting at bridging this gap.

Despite its advantages, the development of robust and scalable solvers for SBM, particularly geometric multigrid
preconditioners~\cite{hackbusch1985multi}, has remained an open problem. This paper addresses this gap by introducing a
geometric multigrid preconditioner specifically designed for SBM. To enable effective smoothing within the multigrid
procedure, in particular to address the non-symmetric and potentially indefinite local structure of cell matrices
adjacent to the surrogate boundary, we introduce the Discontinuous Galerkin (DG)
method~\cite{zienkiewicz2003discontinuous} for SBM discretizations, which, to the best of author's knowledge, has also
not been previously investigated. It is important to note that the choice of DG is motivated by the flexibility it
offers in the construction of the multigrid preconditioner, particularly in designing effective smoothers, rather than
as a general preference over continuous Galerkin (CG) methods. This allows us to leverage the inherent discontinuity of
the DG method to design a cell-wise smoother that effectively addresses the local perturbations introduced by the SBM
discretization. Consequently, the design of a geometric multigrid solver tailored for discontinuous Galerkin SBM
discretizations is the primary contribution of this work, serving as a stepping stone towards efficient solvers for a
wider range of SBM formulations.

In SBM, the \emph{surrogate domain} $\tilde{\Omega}$ is typically constructed as a union of cells from a fixed
background mesh that are deemed \emph{active} (e.g., entirely inside or significantly intersecting the true domain
$\Omega$), and its boundary $\tilde{\Gamma}$ does not conform to the domain boundary $\Gamma$. Boundary conditions are
transferred from the true to the surrogate boundary, typically via Taylor expansions or more general extension
operators~\cite{zorrilla2024shifted}, and enforced in a Nitsche-like manner~\cite{nitsche1971variationsprinzip}. The
Shifted Boundary Method has evolved from its first formulation~\cite{main2018shifted, main2018shiftedVol2}, which used
cells strictly within the considered domain, to a recent approach~\cite{yang2024optimal} that often include intersected
cells based on a volume fraction threshold. It has been extended to high-order discretizations \cite{atallah2022high},
various physical problems including Stokes flow~\cite{atallah2020analysis}, solid mechanics~\cite{atallah2021shifted,
    atallah2024nonlinear}, and significantly, to problems with embedded interfaces~\cite{li2020shifted, xu2024weighted}.
These latter works demonstrate the capability of SBM to handle discontinuities across internal boundaries by
appropriately modifying the formulation to impose jump conditions, extending the method's applicability to multiphysics
and multi-material problems. Other innovations include penalty-free variants~\cite{collins2023penalty} and integration
with level set methods~\cite{kuzmin2022unfitted, xue2021new}.

While this avoids the complexities of generating body-fitted meshes, the primary geometric task in SBM shifts to
accurately determining the relationship between points on the surrogate boundary $\tilde{\Gamma}$ and the true boundary
$\Gamma$. The method inherently allows for the use of arbitrarily complex geometries, and crucially, avoids the need to
compute integrals over the arbitrarily shaped integration domains that arise from cell-boundary intersections. While
determining the active mesh requires computing volume fractions of cut cells, this is a one-time geometric
preprocessing step. In contrast, methods like CutFEM require specialized quadrature for all terms in the bilinear form
on every cut cell. This makes SBM significantly more efficient in terms of computational throughput, especially in
matrix-free implementations where the overhead of cut-cell quadrature would be incurred in every operator evaluation.
The SBM utilizes closest-point projection algorithms to find for each point on $\tilde{\Gamma}$ a corresponding point
on $\Gamma$, but other strategies could be used, such as level sets. Level set methods, which represent the domain
boundary as the zero level set of a function, are often employed in unfitted methods like SBM to facilitate operations
such as closest point projection~\cite{kuzmin2022unfitted, xue2021new}. Special treatment of domains with corners was
analyzed in~\cite{atallah2021analysis}.

However, the geometric flexibility of SBM comes at the cost of new computational challenges. The extrapolation of
boundary conditions and the resulting modifications to the variational formulation often lead to linear systems with
non-symmetric and potentially indefinite properties, particularly for higher-order polynomial approximations, making
their efficient solution non-trivial.

Multigrid methods~\cite{brandt1977multi, hackbusch1985multi} are renowned for their potential to solve large systems
arising from partial differential equations, often offering optimal or near-optimal complexity. However, their
application as preconditioners for SBM remains largely unexplored. While Algebraic Multigrid (AMG) has been applied to
SBM for discretizations using continuous linear elements~\cite{atallah2020second}, its efficiency for higher-order
methods is not established. In fact, AMG generally struggles with high-order finite element discretizations even for
body-fitted meshes, and for unfitted methods like SBM, it is even less effective for $p > 1$. Furthermore, AMG may not
fully leverage the geometric information in structured background meshes. Geometric multigrid methods, in contrast,
explicitly use the hierarchy of meshes and can be very effective, especially when combined with structured or
semi-structured background grids often employed in SBM. A critical component of geometric multigrid is the smoother,
which damps high-frequency error components on each grid level. Standard smoothers like Jacobi or Gauss-Seidel,
however, may struggle with the non-standard local properties of SBM systems occurring at the surrogate boundary. This
motivates the exploration of alternative discretization and smoothing strategies.

Discontinuous Galerkin methods~\cite{reed1973triangular, zienkiewicz2003discontinuous,arnold1982interior,
    cockburn2000development} employ piecewise discontinuous polynomials as basis functions. In interior penalty
formulations~\cite{babuvska1973finite, arnold1982interior}, continuity is enforced weakly by penalizing jumps across
element interfaces. Alternatively, non-symmetric variants exist~\cite{baumann1999discontinuous,
    oden1998discontinuoushpfinite}, which generally exhibit reduced accuracy compared to symmetric
formulations~\cite{zienkiewicz2003discontinuous}. While DG methods typically involve a higher number of degrees of
freedom compared to their continuous Galerkin counterparts, they offer considerable flexibility in return, especially
for implementing Schwarz-like smoothers or $hp$-adaptivity.

The synergy between multigrid techniques and DG discretizations has been a fruitful area of research, leading to
powerful solvers. The inherent block structure of DG methods lends itself naturally to the design of effective
smoothers, a cornerstone of multigrid efficiency. Seminal work, such as the overview in~\cite{antonietti2015multigrid},
has established $hp$-independent convergence rates, a highly desirable property, though the number of smoothing
iterations can sometimes scale with the polynomial degree $p$. Particularly for DG discretizations on structured
Cartesian meshes, specialized smoothers like cell-wise or vertex-patch variants have proven highly effective. When
dealing with separable operators like the Laplacian, the local problems arising in these smoothers can often be solved
with impressive efficiency via fast diagonalization techniques~\cite{witte2021fast,wichrowski2025smoothers,
    wichrowski2025matrix}. This can translate into low iteration counts and robust $hp$-independent convergence.
Furthermore, the adaptability of DG methods to modern hardware architectures is underscored by successful GPU
implementations of multilevel solvers, showcasing their potential for high-performance
computing~\cite{cui2024multilevel}.

The cell-wise and vertex-patch smoothers, mentioned earlier as effective for DG discretizations, can be formally
understood and analyzed within the framework of subspace correction methods. This perspective is particularly valuable
for designing smoothers capable of addressing the local complexities introduced by SBM, such as the non-standard
boundary conditions that can conflict with traditional smoothers. The theory of subspace correction
methods~\cite{xu2001method, bramble1991analysis} provides a powerful framework for analyzing and constructing such
smoothers. While we do not provide a theoretical analysis of the convergence of our multigrid preconditioner, we use
insights from subspace correction theory to guide the design of our smoother. We adopt this perspective and design our
cell-wise smoother as a successive subspace correction method, where each subspace corresponds to the degrees of
freedom on a single cell. The key aspect of this design is the choice of subspace, which allows for effective
accounting for the shifted boundary condition. For practical efficiency, especially with high-order DG methods, the
implementation of these smoothers and the overall multigrid algorithm can significantly benefit from matrix-free
techniques.

Matrix-free implementations significantly enhance the efficiency of DG methods, especially for high-order
discretizations and large problems~\cite{kronbichler2019fast, witte2021fast}. In fact, assembling the system matrix, as
typically required by Algebraic Multigrid (AMG), becomes prohibitively expensive in terms of memory and computational
time for high polynomial degrees. Geometric multigrid methods are particularly well-suited for matrix-free settings, as
their components can often be implemented efficiently without explicit matrix formation~\cite{kronbichler2019fast}.
Recent work~\cite{wichrowski2025MFSBM} has demonstrated that matrix-free evaluation of SBM operators is feasible and
efficient. While this work uses matrix-based implementations to demonstrate the multigrid approach for DG-SBM, it lays
the groundwork for future high-performance matrix-free solvers that avoid matrix assembly altogether.

In the broader context of unfitted methods, CutFEM~\cite{burman2015cutfem} is a prominent alternative with theoretical
results for preconditioners already developed. It discretizes directly on the physical domain by cutting background
cells, requiring specialized quadrature for handling these intersected cells. Due to inherent difficulties with small
cuts, proper stabilization seems to be an unavoidable part of CutFEM. The so-called ghost
penalty~\cite{burman2010ghost, wichrowski2025matrix} solves the issue of ill-conditioning but may require additional
care to avoid locking~\cite{badia2022linking,bergbauer2024high,burman2022design}. CutFEM has been applied to various
problems, including Stokes~\cite{burman2014fictitious}, elasticity~\cite{hansbo2017cut}, or two-phase
flows~\cite{claus2019cutfem}. Furthermore, Discontinuous Galerkin methods have also been combined with
CutFEM~\cite{gurkan2019stabilized, bergbauer2024high}. While CutFEM ensures robust conditioning via geometry-adapted
quadrature, DG-SBM retains the efficiency comparable to standard tensor-product quadrature, allowing for faster
operator evaluation, particularly in matrix-free implementations~\cite{wichrowski2025MFSBM}.

Concerning the preconditioning, CutFEM seems to pose challenges. Results providing optimal
preconditioners~\cite{gross2023analysis, gross2021optimal} have been developed. Although these methods are shown to be
mesh-independent, iteration counts can be high. In~\cite{bergbauer2024high} DG-CutFEM was considered, and a multigrid
preconditioner based on cell-wise Additive Schwarz smoother was used. Although the paper mostly focuses on~matrix-free
implementation, the preconditioner seems promising. However, the smoothing step requires a rather high number of
matrix-vector products.

This paper addresses the computational challenges of solving linear systems arising from SBM discretizations. Our main
contribution is a geometric multigrid preconditioner specifically designed for Discontinuous Galerkin SBM formulations
of the Poisson equation. We develop a cell-wise successive subspace correction smoother that aligns with the DG
framework and the localized nature of SBM. Numerical experiments demonstrate the preconditioner's effectiveness under
mesh refinement, showing strong performance for linear ($p=1$) and quadratic ($p=2$) elements, while identifying
challenges for cubic ($p=3$) elements, particularly in 3D. We also compare our results with an algebraic multigrid
(AMG) preconditioner applied to a standard continuous Galerkin SBM formulation, highlighting the limitations of AMG for
higher-order SBM discretizations. The implementation is built on the \texttt{deal.II} finite element
library~\cite{dealii2019design, dealii2024}, which provides comprehensive tools for finite element methods and
multigrid. The closest point search was implemented specifically for this work.

The remainder of this paper is organized as follows. Section~\ref{sec:method} details the SBM formulation and its DG
discretization. Section~\ref{sec:mg_preconditioner} describes the components of the multigrid preconditioner, with a
focus on the smoother design. Section~\ref{sec:implementation} discusses implementation aspects, including geometry
handling and the closest point projection algorithm. Numerical results evaluating the performance of the proposed
method are presented in Section~\ref{sec:numerical_results}. Finally, Section~\ref{sec:conclusion} offers concluding
remarks and outlines potential directions for future research.

\section{Method formulation}
\label{sec:method}
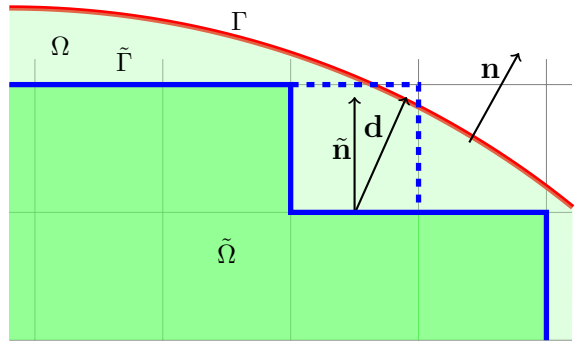
\begin{figure}[!ht]
    \centering
    \begin{tikzpicture}[scale=1.7]
        \draw[step=1cm, gray, very thin] (0.8,0) grid (5.2,2.2);

        \foreach \x/\y in {2/1, 1/0, 1/1, 2/0, 3/0, 4/0}
            {
                \fill[green, opacity=0.5] (\x,\y) rectangle (\x+1,\y+1);
            }
        \foreach \x/\y in {1/0 }{
                \fill[green, opacity=0.5] (\x,\y) rectangle (\x-0.2,\y+2);
            }

        \draw[line width=2pt, red, domain=51:90, samples=100, variable=\t]
        plot ({0.8 + 7.0*cos(\t)}, {-4.4 + 7.0*sin(\t)});
        \node at (2.6,2.5) {$\Gamma$};

        \fill[green!30, opacity=0.4] plot[domain=51:90, samples=100, variable=\t] ({0.8 + 7.0*cos(\t)}, {-4.4 +
                7.0*sin(\t)}) -- (0.8,0) -- (5.2,0) -- cycle;

        \pgfmathsetmacro{\myangle}{66}
        \draw[<-, thick] ({0.8 + 6.9*cos(\myangle) + 0.3}, {-4.4 + 6.9*sin(\myangle)}) -- ({0.8 + 5.9*cos(\myangle)
                +
                0.3}, {-4.4 +
                5.9*sin(\myangle)});
        \pgfmathsetmacro{\labelangle}{\myangle - 6} 
        \node at ({0.8 + 6.9*cos(\labelangle) - 0.6}, {-4.4 + 6.9*sin(\labelangle)+0.1}) {\large $\mathbf{d}$};

        \draw[->, thick] ({0.8 + 6.8*cos(\myangle-5) + 0.3}, {-4.4 + 6.8*sin(\myangle-5)}) -- ({0.8 +
                7.6*cos(\myangle-5)
                +
                0.3}, {-4.4 +
                7.6*sin(\myangle-5)});
        \node at ({0.8 + 7.5*cos(\labelangle-5) - 0.54}, {-4.4 + 7.5*sin(\labelangle-5)+0.36}) {\large $\mathbf{n}$};

        \draw[<-, thick]  ({0.8 + 5.9*cos(\myangle) +
                0.3}, {-4.4 + 6.9*sin(\myangle)}) -- ({0.8 + 5.9*cos(\myangle) +
                0.3}, {-4.4 +
                5.9*sin(\myangle)});
        \node at (3.4, 1.5) {\large $\tilde{\mathbf{n}}$}; 

        \draw[line width=2pt, blue] (0.8,2) -- (3,2) -- (3,1) -- (4,1) -- (5,1)-- (5,0);

        \draw[line width=2pt, blue, dashed] (3,2) -- (4,2) -- (4,2) -- (4,1) ;
        \node at (1.7,2.2) {$\tilde{\Gamma}$};

        \node at (2.5,0.7) {$\tilde{\Omega}$};

        \node at (1.2,2.3) {${\Omega}$};

    \end{tikzpicture}
    \caption{Schematic illustrating the background mesh, interior cells (green), the surrogate boundary
        $\tilde{\Gamma}$ (thick blue line) along the upper boundary of the interior cells, and the true boundary
        $\partial\Omega$ (red).}
    \label{fig:sbm_setup}
\end{figure}

We consider the Poisson problem as a model problem:
\begin{subequations}
    \begin{align}
        -\Delta u & = f \quad \text{in } \Omega,  \label{eq:poisson}      \\
        u         & = g \quad \text{on } \Gamma,  \label{eq:dirichlet_bc}
    \end{align}
\end{subequations}
where $\Omega \subset \mathbb{R}^d$ ($d=2,3$) is a domain with boundary $\Gamma =\partial\Omega$ as depicted in
Figure~\ref{fig:sbm_setup}, $f$ is a given source term, and $g$ is the prescribed Dirichlet boundary condition.

We restrict our attention to Dirichlet boundary conditions deliberately. The contribution of this work is a geometric
multigrid \emph{preconditioner} for SBM systems, and the algebraic difficulty it must overcome---the loss of symmetry
and the possible indefiniteness introduced by the boundary \emph{extrapolation}---is already fully present in the
Dirichlet case. Neumann (and Robin) conditions, by contrast, raise a distinct \emph{discretization} question: a naive
shifted enforcement of the flux condition reduces the convergence order, and recovering optimality requires a dedicated
treatment, as recently analyzed and resolved in~\cite{COLLINS2026118793}. Incorporating an under-resolved Neumann
discretization would confound discretization error with the solver behavior we aim to isolate, while adopting an
optimal-order Neumann formulation is largely orthogonal to the preconditioner design. Since the cell-wise smoother and
the transfer operators developed below act on the algebraic system and are agnostic to the type of weakly imposed
boundary condition, the proposed multigrid framework carries over once such a formulation is in place; we therefore
defer the Neumann and mixed cases to future work.

To solve this problem numerically, we first formulate it in a weak sense. We seek a solution $u$ in an appropriate
function space, $H^1(\Omega)$, which consists of functions that are square-integrable and whose first derivatives are
also square-integrable. Multiplying the equation by a test function $v \in H^1(\Omega)$ and integrating over $\Omega$,
we obtain:
\begin{equation}
    \int_\Omega \nabla u \cdot \nabla v \, dx - \int_{\Gamma} (\nabla u \cdot \mathbf{n}) \, v \, ds =
    \int_\Omega f \, v \, dx.
\end{equation}

We next introduce a mesh \( \mathcal{T}_h \) consisting of quadrilateral (2D) or hexahedral (3D) elements of size $h$,
and define a finite element space $\mathbb{V}_h\subset H^1(\Omega)$ using Lagrange polynomial elements of degree \( p
\). In classical finite element methods, the mesh conforms to the boundary (unlike the background mesh approach
illustrated in Figure~\ref{fig:sbm_setup}), and test functions $v$ typically vanish on $\Gamma $ to strongly enforce
the Dirichlet condition $u=g$. When function spaces do not necessarily satisfy essential boundary conditions strongly,
boundary conditions must be enforced weakly through additional integral terms on $\Gamma$.

Nitsche's method provides a way to weakly impose Dirichlet boundary conditions within a variational formulation without
requiring the function space to satisfy the boundary conditions. It modifies the bilinear form by adding terms on the
boundary $\Gamma$. The standard Nitsche formulation for the Poisson problem with Dirichlet boundary conditions $u=g$ on
$\Gamma$ seeks $u_h \in \mathbb{V}_h$ such that for all $v_h \in \mathbb{V}_h$:
\begin{equation}
    \begin{split}
        \int_\Omega \nabla u_h \cdot \nabla v_h \, dx
        \;\; - \int_{\Gamma} (\nabla u_h \cdot \mathbf{n}) v_h \, ds
        - & \alpha\int_{\Gamma} (\nabla v_h \cdot \mathbf{n}) u_h \, ds
        + \int_{\Gamma} \sigma_{\Gamma} u_h v_h \, ds =                 \\
          & = \int_\Omega f v_h \, dx
        - \alpha \int_{\Gamma} (\nabla v_h \cdot \mathbf{n}) g \, ds
        + \int_{\Gamma} \sigma_\Gamma g v_h \, ds.
    \end{split}
\end{equation}
Here, $\sigma_\Gamma$ is a penalty parameter, typically chosen as $\sigma_\Gamma = \mathcal{O}(p^2 h^{-1})$ for mesh
size $h$, and
$\mathbf{n}$ is the outward unit normal vector to $\Gamma$.  The choice of parameter $\alpha=1$ leads to a symmetric
formulation, while $\alpha=-1$ results in a non-symmetric formulation in which the penalty term can be
skipped~\cite{baumann1999discontinuous}.

\subsection{Shifted Boundary Method}
In many applications, the domain $\Omega$ may have a complex geometry, making the generation of body-fitted meshes
challenging. Non-body-fitted (unfitted) methods address this challenge by employing a background mesh $\mathcal{T}_h$
that does not conform to the boundary of $\Omega$. The domain $\Omega$ is embedded within this background mesh, and
cells of $\mathcal{T}_h$ are classified as active based on their intersection with $\Omega$. The computational domain
$\tilde{\Omega}$ is defined as the union of these active cells. In the original SBM formulation~\cite{main2018shifted},
only cells strictly contained in $\Omega$ were included, leading to the surrogate domain boundary depicted by the solid
blue line in Figure~\ref{fig:sbm_setup}. Consequently, the surrogate domain $\tilde{\Omega}$ is generally a subset of
$\Omega$, and its boundary $\tilde{\Gamma}$ does not coincide with the true boundary $\Gamma$.

Later extensions include intersected cells~\cite{yang2024optimal}. For each intersected cell $K$, we define the
\emph{volume fraction}
\begin{equation}
    \kappa(K) = \frac{|K \cap \Omega|}{|K|} \in [0,1],
\end{equation}
i.e.\ the fraction of the cell measure (area in 2D, volume in 3D) that lies inside the true domain $\Omega$. A cell is
classified as active if $\kappa(K) > \lambda$, where the threshold $\lambda \in [0,1]$ controls how much of a cell must
be inside $\Omega$ for it to be retained. The two limits are instructive: $\lambda = 1$ recovers the original
formulation that keeps only cells strictly inside $\Omega$, while $\lambda \to 0$ admits every cell touched by the
domain. Intermediate values trade off how closely the surrogate boundary $\tilde{\Gamma}$ approximates $\Gamma$ against
the conditioning of the resulting system: a small $\lambda$ reduces the distance between the true and surrogate
boundaries by including sliver cells with little volume inside $\Omega$, but such cells can lead to ill-conditioning. In
Figure \ref{fig:sbm_setup}, lowering $\lambda$ corresponds to including the additional cells indicated by the dashed
blue line. We note that $\kappa(K)$ need not be evaluated exactly: in this work we use accurate non-matching quadrature
(Section~\ref{sec:implementation}), but a cheaper surrogate—such as the fraction of cell nodes at which the level set
has the interior sign—is sufficient to classify cells, since the threshold $\lambda$ only enters through the inequality
$\kappa(K) > \lambda$.

To impose boundary conditions on $\tilde{\Gamma}$, the Dirichlet condition prescribed on the true boundary $\Gamma$ is
extrapolated to the surrogate boundary via an \emph{extension operator} $\mathcal{E}$, which maps functions defined on
$\tilde{\Omega}$ to the surrogate boundary $\tilde{\Gamma}$. For each point $\tilde{\mathbf{x}} \in \tilde{\Gamma}$ on
the surrogate boundary, let $\mathbf{x} \in \Gamma$ be its closest-point projection onto the true boundary, and let
$\mathbf{d} = \mathbf{x} - \tilde{\mathbf{x}}$ denote the shift vector. We assume that the Dirichlet data $g$ on
$\Gamma$ is the restriction of a smooth function $u^\star$ defined in a neighborhood of $\tilde{\Gamma}$, so that
$u^\star(\mathbf{x}) = g(\mathbf{x})$. A first-order Taylor expansion of $u^\star$ about $\tilde{\mathbf{x}}$ gives
\begin{equation}
    u^\star(\mathbf{x}) = u^\star(\tilde{\mathbf{x}}) + \mathbf{d} \cdot \nabla u^\star(\tilde{\mathbf{x}}) + R,
\end{equation}
where $R$ collects the higher-order remainder. Solving for the value at the surrogate boundary and dropping the
remainder defines the extension operator $\mathcal{E}$:
\begin{equation}
    \mathcal{E}u^\star(\tilde{\mathbf{x}}) := u^\star(\mathbf{x}) - \mathbf{d} \cdot \nabla u^\star(\tilde{\mathbf{x}})
    = g(\mathbf{x}) - \mathbf{d} \cdot \nabla u^\star(\tilde{\mathbf{x}}).
    \label{eq:extension_operator}
\end{equation}
In the discrete scheme, $u^\star$ is taken to be the discrete solution $u_h$ in the surrogate domain $\tilde{\Omega}$,
so that $\nabla u^\star(\tilde{\mathbf{x}}) = \nabla u_h(\tilde{\mathbf{x}})$ is evaluated directly from the polynomial
representation of $u_h$, avoiding the computation of higher-order derivatives. Substituting $\mathcal{E}u_h$ into the
weak formulation yields a variational problem in which the extrapolated boundary condition on $\tilde{\Gamma}$ is
enforced in a Nitsche-like manner.

The SBM weak formulation seeks $u_h \in \mathbb{V}_h$ such that for all $v_h \in \mathbb{V}_h$,
\begin{equation}
    \begin{split}
        \int_{\tilde{\Omega}} \nabla u_h \cdot \nabla v_h \, dx
        - \int_{\tilde{\Gamma}} (\nabla u_h \cdot \tilde{\mathbf{n}}) v_h \, ds
        - \alpha \int_{\tilde{\Gamma}} (\nabla v_h \cdot \tilde{\mathbf{n}}) \; \mathcal{E}u_h \, ds
        + \int_{\tilde{\Gamma}} \sigma_\Gamma \; \mathcal{E} u_h  \; v_h \, ds = \\
        = \int_{\tilde{\Omega}} f v_h \, dx
        - \alpha \int_{\tilde{\Gamma}} (\nabla v_h \cdot \tilde{\mathbf{n}}) \, g(\tilde{\mathbf{x}} + \mathbf{d}) \, ds
        + \int_{\tilde{\Gamma}} \sigma_\Gamma \, g(\tilde{\mathbf{x}} + \mathbf{d}) \, v_h \, ds.
    \end{split}
    \label{eq:sbm_weak}
\end{equation}
where $\tilde{\mathbf{n}}$ is the outward unit normal to $\tilde{\Gamma}$. On the right-hand side, the Dirichlet datum
$g$ is defined only on the true boundary $\Gamma$; at a point $\tilde{\mathbf{x}} \in \tilde{\Gamma}$ it is therefore
evaluated at the closest-point projection $\mathbf{x} = \tilde{\mathbf{x}} + \mathbf{d} \in \Gamma$, which we indicate
explicitly by writing $g(\tilde{\mathbf{x}} + \mathbf{d})$. This data term replaces the extrapolated trial term
$\mathcal{E}u_h$ appearing on the left-hand side: by construction, the extension operator reproduces the true boundary
data for the exact solution, $\mathcal{E}u^\star(\tilde{\mathbf{x}}) = u^\star(\mathbf{x}) = g(\tilde{\mathbf{x}} +
    \mathbf{d})$. Note that, unlike $\mathcal{E}u_h$ in~\eqref{eq:extension_operator}, the data term contains no gradient
contribution, since the extrapolation is exact for the prescribed boundary values.
Note that the resulting form is no longer symmetric due to the presence of the extrapolated
boundary condition. We will refer to the formulation with $\alpha=1$ as the \emph{quasi-symmetric} formulation.

Moreover, a simple 1D experiment using the domain $\Omega=[0, \xi]$ with a surrogate domain $\tilde{\Omega} = [0, 1]$
and a single cell shows that for Lagrange polynomials of degree $p > 1$, the stiffness matrix can have negative entries
on the diagonal. As a result, the system matrix may lose diagonal dominance, which can adversely affect the stability
and convergence of iterative solvers.

\begin{figure}[!ht]
    \centering
    \begin{tikzpicture}
        \begin{axis}[
                xlabel=$\xi$,
                ylabel=Imaginary part of eigenvalue,
                legend style={at={(0.9,0.82)},anchor=east},
                grid=major,
                xmin=0.4, xmax=2,
                height=8cm, 
            ]

            \addplot[red, thick, mark=*, mark repeat=14] table[x index=0, y index=1]
                {results/standard_1d/imagin_eigs1_p1.csv};
            \addlegendentry{$p=1$}
            \addplot[darkgreen, thick, mark=triangle*, mark repeat=16] table[x index=0, y index=1]
                {results/standard_1d/imagin_eigs1_p2.csv};
            \addlegendentry{$p=2$}
            \addplot[blue, thick, mark=+, mark repeat=15] table[x index=0, y index=1]
                {results/standard_1d/imagin_eigs1_p3.csv};
            \addlegendentry{$p=3$}

            \addplot[red, thick, dashed, mark=*, mark repeat=13] table[x index=0, y index=1]
                {results/penalty_free/imagin_eigs1_p1.csv};
            \addplot[darkgreen, thick, dashed, mark=triangle*, mark repeat=17] table[x index=0, y index=1]
                {results/penalty_free/imagin_eigs1_p2.csv};
            \addplot[blue, thick, dashed, mark=+, mark repeat=15] table[x index=0, y index=1]
                {results/penalty_free/imagin_eigs1_p3.csv};

            \draw[gray, thick] (axis cs:1, \pgfkeysvalueof{/pgfplots/ymin}) -- (axis cs:1,
            \pgfkeysvalueof{/pgfplots/ymax});
        \end{axis}
    \end{tikzpicture}
    \caption{Maximal imaginary part of eigenvalues for the 1D SBM stiffness matrix versus true domain size $\xi$. The
        surrogate domain is $[0,1]$ (single cell); $\xi=1$ implies coincident boundaries. Solid lines: quasi-symmetric
        penalty ($\alpha=1$, $\sigma_\Gamma=5$); dashed lines: non-symmetric penalty ($\alpha=-1$,
        $\sigma_\Gamma=0$). }
    \label{fig:1d_eigenvalues}
\end{figure}
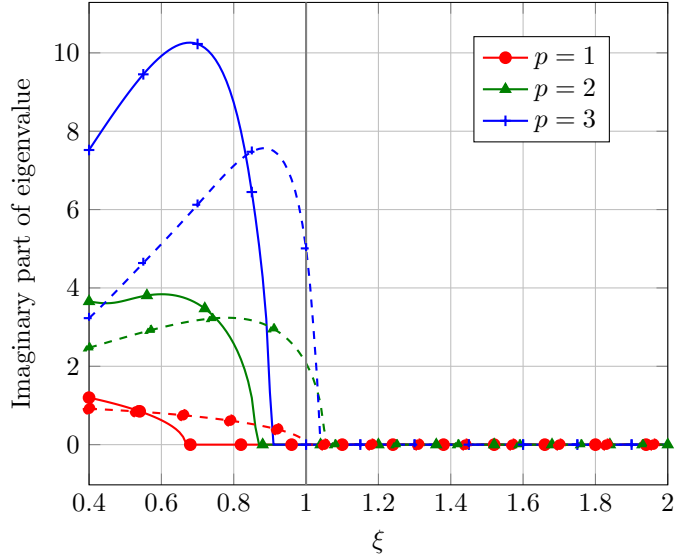

Notably, a penalty-free formulation ($\sigma_\Gamma=0$) can be achieved by setting $\alpha=-1$, which alters the sign
of the extrapolated term in the weak formulation~\cite{collins2023penalty}. This approach also results in a positive
diagonal for the stiffness matrix, and both formulations yield optimal convergence rates. However, symmetry of the
remaining terms in the weak formulation is not recovered. Consequently, the 1D system matrix for this problem can
exhibit complex eigenvalues when the true boundary $\xi$ is less than the surrogate boundary (i.e., $\xi < 1$). Figure
\ref{fig:1d_eigenvalues}~illustrates this by plotting the maximum imaginary part of the eigenvalues for this 1D problem
across different polynomial degrees $p$.

Figure~\ref{fig:1d_eigenvalues} reveals that complex eigenvalues with non-zero imaginary parts indeed arise for
negative shifts (when the true domain boundary $\xi$ is less than the surrogate domain boundary, $\xi < 1$). The
penalty-free formulation (dashed lines) exhibits eigenvalues with larger imaginary parts compared to~the
quasi-symmetric penalty formulation (solid lines) in this regime, while for positive shifts ($\xi > 1$), the
eigenvalues remain real for both formulations. The presence of complex eigenvalues is problematic for multigrid
methods, as they can adversely affect smoothing properties and overall convergence of the preconditioner. Due to the
tensor product structure of the Laplace operator, these 1D observations suggest that similar issues with complex
eigenvalues may manifest in higher dimensions (2D and 3D), potentially compromising the effectiveness of standard
smoothers and degrading multigrid performance.

Interestingly, for higher-order polynomials ($p>1$), the penalty-free formulation can yield eigenvalues with nonzero
imaginary parts even when the surrogate boundary coincides with the true boundary ($\xi=1$). In contrast, the
quasi-symmetric formulation is generally less prone to such imaginary eigenvalues. For this reason, we choose the
quasi-symmetric formulation as the default in our work.

These observed spectral properties, particularly the emergence of complex eigenvalues, stem from the SBM's treatment of
boundary conditions on the surrogate boundary $\tilde{\Gamma}$. The effects of this treatment are localized to the
cells adjacent to $\tilde{\Gamma}$ where the extrapolation and non-standard penalty terms are applied. Intuitively, one
might expect that the \emph{Laplacian-like} properties are perturbed only within these boundary-adjacent cells. This
localization suggests that it should be possible to resolve these local perturbations within the scope of individual
cells. Discontinuous Galerkin (DG) methods offer a natural framework for such a cell-wise treatment, as they inherently
allow for different solution behaviors and formulations on each element and provide mechanisms to couple them
appropriately.

\subsection{Discontinuous Galerkin Discretization}

We define the DG finite element space $\DG$ as piecewise polynomial functions of degree $p$ on the mesh
$\mathcal{T}_h$, without continuity requirements across element boundaries. The DG weak formulation extends the
standard variational formulation by incorporating terms that account for discontinuities across interior faces
$\mathcal{F}_h^{\text{int}}$.

For a scalar function $w$, we define the jump operator $[w]$ and the average operator $\{w\}$ across an interior face
$F \in \mathcal{F}_h^{\text{int}}$ as:
\begin{align}
    [w]   & = w^+ \mathbf{n}^+ + w^- \mathbf{n}^-, \\
    \{w\} & = \frac{1}{2}(w^+ + w^-),
\end{align}
where $w^+$ and $w^-$ denote the values of $w$ on the two elements sharing the face $F$, and $\mathbf{n}^+$ and
$\mathbf{n}^-$ are the outward unit normal vectors to $F$ on these elements. In the symmetric interior penalty
discontinuous Galerkin method, the bilinear form is modified to include these jump and average terms as follows:
\begin{equation}
    \begin{split}
        \int_{\tilde{\Omega}} \nabla u_h \cdot \nabla v_h \, dx \longrightarrow   \int_{\mathcal{T}_h}
        \nabla
        u_h \cdot \nabla v_h \, dx
         & - \sum_{F \in \mathcal{F}_h^{\text{int}}} \int_F \big( \{ \nabla u_h \} \cdot [v_h] + \{ \nabla v_h \} \cdot
        [u_h] \big) \, ds                                                                                                 \\
         & + \sum_{F \in \mathcal{F}_h^{\text{int}}} \int_F \sigma_F [u_h]\cdot [v_h] \, ds = a_{\tilde{\Omega}}(u_h, v_h
        ).
    \end{split}
    \label{eq:dg_bilinear}
\end{equation}
The first term represents the standard DG contribution, while the second and third terms account for the
discontinuities across the interior faces. The penalty term involving the parameter $\sigma_F$ ensures stability of the
method.

By putting the form $a_{\tilde{\Omega}} (u_h, v_h)$ from the equation above together with the Nitsche-like term, we
obtain the complete bilinear form for the DG-SBM discretization of the Poisson problem. The problem seeks $u_h \in \DG$
such that for all $v_h \in \DG$:
\begin{equation}
    \begin{split}
        a_h(u_h, v_h) =    a_{\tilde{\Omega}} (u_h, v_h) - \int_{\tilde{\Gamma}} (\nabla u_h \cdot \tilde{\mathbf{n}})
        v_h \, ds
        - & \alpha\int_{\tilde{\Gamma}} (\nabla v_h \cdot \tilde{\mathbf{n}}) \mathcal{E}u_h \, ds
        + \int_{\tilde{\Gamma}} \sigma_\Gamma \mathcal{E} u_h v_h \, ds =                          \\
          & = \int_{\tilde{\Omega}} f v_h \, dx
        - \alpha \int_{\tilde{\Gamma}} (\nabla v_h \cdot \tilde{\mathbf{n}}) \, g(\tilde{\mathbf{x}} + \mathbf{d}) \, ds
        + \int_{\tilde{\Gamma}} \sigma_\Gamma \, g(\tilde{\mathbf{x}} + \mathbf{d})  v_h \, ds,
    \end{split}
\end{equation}
where, as in~\eqref{eq:sbm_weak}, the boundary datum is evaluated at the projected point $\tilde{\mathbf{x}} +
    \mathbf{d} \in \Gamma$.
By choosing a basis $\{\phi_i\}_{i=1}^N$ for $\DG$ in a standard form, we arrive at a linear system $A_h u_h = f_h$,
where $(A_h)_{ij} =
    a_h(\phi_i, \phi_j)$, $f_h$ is the right-hand side vector, and $u_h$ is the vector of coefficients. Since~the
basis defines a one-to-one mapping between real vectors and functions from $\DG$, we will not distinguish between these
two. While alternative orthogonal bases such as Legendre or Chebyshev polynomials are often recommended to improve the numerical stability and conditioning of high-order DG formulations~\cite{cockburn1998local, hesthaven2008nodal, jaskowiec2018very, jaskowiec2017application, jaskowiec2024penalty}, our numerical experiments indicate that replacing the nodal Lagrange basis with Legendre polynomials yields almost no changes in the overall multigrid iteration counts. This is likely because we consider relatively low polynomial degrees (up to $p=3$), where the conditioning differences between the bases are negligible. For higher orders, the overall solver performance is currently dominated by other structural challenges of the SBM formulation rather than the choice of local basis. Thus, we maintain the use of Lagrange polynomials.

\section{Multigrid Preconditioner}
\label{sec:mg_preconditioner}
The DG-SBM formulation described above leads to a large, sparse linear system. Due to
the lack of symmetry and possible indefiniteness of the resulting matrix, we employ a~Krylov subspace method,
specifically GMRES, for its solution.  The convergence of GMRES is sensitive to the condition number of the system
matrix. For symmetric positive definite discretizations of the Poisson problem, this condition number scales as
$\mathcal{O}(p^{4}\, h^{-2})$ for the interior penalty DG method~\cite{antonietti2015multigrid}, deteriorating both as
the mesh is refined (i.e.\ as the number of elements grows) and as the polynomial degree $p$ increases. The SBM
extrapolation terms further degrade this conditioning near the surrogate boundary.
Multigrid methods~\cite{hackbusch1985multi} rely on a hierarchy of meshes with varying resolution levels. Our Cartesian
background mesh
naturally facilitates the construction of such a nested hierarchy. We define a sequence of
meshes $\{\mathcal{T}_\ell\}_{\ell=0}^L$, where $\ell$ denotes the level and $L$ is the finest level. For each level,
we identify active cells that have at least a fraction $\lambda$ of their measure (area in 2D, volume in 3D) inside the
domain.
On each mesh $\mathcal{T}_\ell$, we define a finite element space $\DGell$ consisting of piecewise polynomials of~
degree $p_\ell$.  Since the set of active cells is determined by their intersection with the domain boundary, these
sets are not necessarily nested between levels, which can affect the efficiency of standard inter-grid transfer
operators. While our construction ensures a nested hierarchy  of finite element
spaces, the lack of a strict geometric relationship between active cells at different levels may reduce the
effectiveness of transfer operators. To address this, we employ a $p$-multigrid strategy combined with $h$-multigrid: we use only linear elements ($p=1$) on
geometrically coarser meshes and gradually increase the polynomial degree from $p=1$ to the target degree $p$ only on the finest mesh. This approach, often referred to as $p$-preconditioning, results in a total number of levels equal to $L+p$.

With this hierarchy of meshes and spaces established, the multigrid method is defined through two crucial components:
the smoothers, which reduce high-frequency error components on each level, and the transfer operators, which move
information between different resolution levels. These components, together with a coarse-grid solver, form the
complete multigrid algorithm. For preconditioning, we use a single multigrid V-cycle.

\subsection{Smoother}
\label{sec:smoother}
We first consider Richardson iterations with preconditioners $P$. Let $A_h$ denote the system matrix arising from the discretization on the current mesh level. Given the current
approximation $u_h$ and the right-hand side $f_h$, a~smoothing step updates
the approximation:
\begin{equation}
    u_h^{k+1} = u_h^k + P(f_h - A_h u_h^k).
\end{equation}
We decompose the finite element space $\DGell$ into $J$ subspaces $\mathbb{V}_i$, each spanned by a subset of the global
basis $\{\phi_j\}$, so that every $v \in \DGell$ can be written as $v = \sum_{i=1}^J v_i$ with $v_i \in \mathbb{V}_i$.
Note that the number of subspaces $J$ is in general different from the dimension $N$ of $\DGell$ (the number of global
basis functions); for the cell-wise choice adopted below, $J$ equals the number of cells. We denote this by the
\emph{subspace sum}
\begin{equation}
    \DGell = \sum_{i=1}^J \mathbb{V}_i,
\end{equation}
which is the span of the union $\bigcup_{i=1}^J \mathbb{V}_i$; we write a sum rather than a set union because a union of
subspaces is in general not itself a subspace, and because the decomposition need not be direct (the subspaces may
overlap). It is convenient to
work in coordinates: identifying $\DGell$ with $\mathbb{R}^N$ through this basis, the subspace $\mathbb{V}_i$ of
dimension $n_i$ corresponds to $\mathbb{R}^{n_i}$, and the embedding is represented by the Boolean selection matrix
$R_i^T \in \mathbb{R}^{N \times n_i}$ that maps a local coefficient vector to a global one by inserting it in the
positions belonging to $\mathbb{V}_i$ and padding with zeros. Its transpose $R_i \in \mathbb{R}^{n_i \times N}$ is the
restriction that picks out those entries. Working with coordinate vectors makes the operations unambiguous and avoids
the identification of $\DGell$ with its dual: $R_i$ and $R_i^T$ act on coefficient vectors, not on functionals. The
local operator $A_i \in \mathbb{R}^{n_i \times n_i}$ is the Galerkin restriction of the global matrix to the subspace,
\begin{equation}
    A_i = R_i A_h R_i^T,
    \label{eq:local_operator}
\end{equation}
i.e.\ the principal submatrix of $A_h$ indexed by the degrees of freedom of $\mathbb{V}_i$; equivalently, $(A_i)_{kl} =
    a_h(\phi_{i,l}, \phi_{i,k})$, where $\{\phi_{i,k}\}$ is the basis of $\mathbb{V}_i$. Then, the additive subspace
correction preconditioner is defined as:
\begin{equation}
    P = \omega \sum_{i=1}^J R_i^T A_i^{-1} R_i,
\end{equation}
where $\omega$ is a relaxation parameter. Alternatively, the successive subspace correction method visits the subspaces
$\mathbb{V}_i$ in a sequential manner. The preconditioner updates the solution by sequentially applying corrections for
each subspace. For each subspace $\mathbb{V}_i$, a correction is computed and applied:
\begin{equation}
    u_h \leftarrow u_h + R_i^T A_i^{-1} R_i (f_h - A_h u_h).
\end{equation}
This process is repeated for all subspaces $i=1, \ldots, J$. If the subspaces are chosen as one-dimensional spaces
spanned by the basis functions, then the preconditioner is equivalent to either Jacobi (additive) or Gauss-Seidel
(successive).

In the context of DG, a natural choice for the subspaces $\mathbb{V}_i$ is the set of functions that are non-zero only
on a single cell. This leads to a cell-wise smoother, where we correct the solution on each cell independently. This
can be interpreted as a block Gauss-Seidel smoother, where each block corresponds to the degrees of freedom associated
with a single cell. With this choice, each subspace is small and of fixed size: $\dim \mathbb{V}_i = {(p_\ell + 1)}^d$,
independent of the number of cells, so that $A_i$ is a small dense block whose factorization is inexpensive and is
reused at every smoothing step. The cost of forming and inverting all blocks therefore scales linearly with the number
of cells, which is what makes the local solves computationally worthwhile. For high polynomial degrees, where the dense
block $A_i$ grows as ${(p_\ell+1)}^{2d}$ and its direct inversion becomes comparatively costly, the local solve
$A_i^{-1}$ could itself be replaced by an inexact local solver, for instance a $p$-multigrid sweep on the
cell~\cite{wichrowski2025local}; we use the exact dense inverse here, but this is a promising direction for improving
the per-cell cost at higher orders.

The local matrix $A_i$ is non-singular for the cell-wise choice. Each diagonal block is the operator $a_h$ restricted
to the polynomials supported on a single cell $K$, which includes the cell's volume term $\int_K \nabla u \cdot \nabla
    v$ together with the interior-penalty face terms that, on the cell's own faces, contribute the symmetric positive
semi-definite jump penalties; for active cells the resulting block is invertible, and in the SBM case the surrogate
boundary contributes additional Nitsche terms that we account for in the local solve. Should a block become nearly
singular—for example on a sliver active cell with a small volume fraction—we stabilize it as discussed in
Section~\ref{sec:implementation}. We emphasize that, with the coordinate definitions above, the preconditioner $P$ is
not an ad hoc heuristic but the standard additive (respectively multiplicative) subspace correction operator for the
decomposition $\DGell = \sum_i \mathbb{V}_i$; the formula $A_i = R_i A_h R_i^T$ in~\eqref{eq:local_operator} is exactly
the Galerkin projection of $A_h$ onto $\mathbb{V}_i$, and coincides with the change-of-basis derivation in which
$R_i^T$ relates the global and local bases.

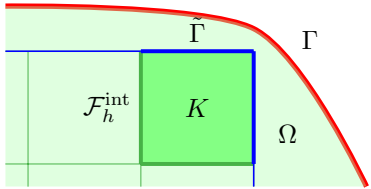
\begin{figure}[!ht]
    \centering
    \begin{tikzpicture}[scale=1.5]

        \draw[step=1cm, darkgreen, very thin] (-1.2,-0.2) grid (1.,1.);
        \draw[fill=green, opacity=0.6] (0,0) rectangle (1,1);

        \draw[line width=1.5pt, darkgreen] (0,0) -- (1,0); 
        \draw[line width=1.5pt, darkgreen] (0,0) -- (0,1); 

        \draw[line width=2pt, red] plot [smooth] coordinates {(2,-0.2) (1,1.2) (-1.2, 1.4)};
        \node at (1.5, 1.1) {$\Gamma$};

        \fill[green!30, opacity=0.4] (-1.2, -0.2) -- (2,-0.2) -- plot [smooth] coordinates {(2,-0.2) (1,1.2) (-1.2,
                1.4)}
        -- cycle;

        \draw[line width=1.5pt, blue] (0,1) -- (1,1); 
        \draw[line width=1.5pt, blue] (1,0) -- (1,1); 

        \draw[line width=0.6pt, blue] (-1.2,1) -- (0,1); 
        \draw[line width=0.6pt, blue] (1,0) -- (1,-0.2); 

        \node[left] at (0,0.5) {$  \mathcal{F}_h^{\text{int}} $};
        \node[above] at (0.5,1) {$ \tilde{\Gamma}$};

        \node[above] at (1.3,0.1) {$ {\Omega}$};
        \node at (0.5,0.5) {$ K$};

    \end{tikzpicture}
    \caption{
        Illustration of a cell $K$ within the DG-SBM framework, highlighting different face types relevant to the cell-wise
        smoother. Interior faces ($\mathcal{F}_h^{\text{int}}$, green) connect to neighboring active cells. For cells
        adjacent
        to the domain boundary, faces can be part of the surrogate boundary $\tilde{\Gamma}$ (blue), which is shifted from
        the
        true boundary $\Gamma$ (red). The cell-wise smoother solves a local problem on $K$, implicitly handling these
        conditions.
    }
    \label{fig:single_cell_boundaries}
\end{figure}

Viewing the smoother as a subspace correction method reveals important insights into the convergence properties. The
correction computed within a subspace $\mathbb{V}_i$ implicitly satisfies certain boundary conditions on the boundary
of that subspace. This is particularly important in the context of the SBM, where the boundary conditions are shifted.
A careless choice of subspace can lead to overconstraining the correction, preventing certain error components from
being effectively smoothed. For instance, Jacobi and Gauss-Seidel smoothers for continuous finite element methods,
which correspond to vertex-patch smoothing with zero correction at patch boundaries, conflict with the shifted boundary
conditions and are not effective, as confirmed in our preliminary numerical experiments. In contrast, the cell-wise
smoother aligns well with the DG method's element-based structure and takes into account the local nature of the SBM
boundary conditions, correcting the solution on each cell independently, as illustrated in Figure
\ref{fig:single_cell_boundaries}.

We adopt this perspective and design our cell-wise smoother as a successive subspace correction method. Specifically,
we implement a cell-wise Symmetric Successive Over-Relaxation (SSOR) smoother, which is a variant of the block
Gauss-Seidel method where each block corresponds to the degrees of freedom on a single cell.

We stress that the locality of the per-cell solves does not preclude effective smoothing. Two mechanisms provide the
necessary inter-cell communication. First, in the multiplicative (Gauss-Seidel/SSOR) variant the cells are visited
sequentially and each correction uses the already-updated values of the previously visited neighbors, so information
propagates across the mesh within a single sweep. Second, and more fundamentally, the DG operator itself couples
neighboring cells through the interior-penalty face terms $\{\nabla u\}$, $[u]$ in~\eqref{eq:dg_bilinear}: although the
local solve is restricted to one cell, the residual that drives it carries the contribution of the neighbors' current
values through these face terms. Cell-wise (block-Jacobi and block-Gauss-Seidel) and vertex-patch smoothers are, in
fact, a well-established and highly effective family for DG discretizations of elliptic
problems~\cite{kronbichler2019fast, witte2021fast, antonietti2015multigrid}, and the purely \emph{additive}
block-Jacobi cell-wise smoother has been used successfully even for the closely related DG-CutFEM
setting~\cite{bergbauer2024high}. Our choice of the multiplicative SSOR variant is motivated by robustness with respect
to the relaxation parameter rather than by a deficiency of the additive form.

\subsection{Transfer Operators}
\label{sec:transfer_operators}
The transfer operators move data between the finite element spaces associated with consecutive multigrid levels. We
define them on the spaces over the \emph{full background mesh}, $\mathbb{W}_\ell^{\text{DG}} \supseteq
    \mathbb{V}_\ell^{\text{DG}}$, rather than on the surrogate-domain spaces $\mathbb{V}_\ell^{\text{DG}}$ directly. This
distinction matters: the background-mesh spaces are nested across an $h$-refinement step (each coarse cell is the union
of its children and the coarse polynomial space is contained in the fine one), whereas the surrogate-domain spaces
$\mathbb{V}_\ell^{\text{DG}}$ are \emph{not} nested. Because SBM discards intersected cells, refinement lets smaller
cells fit between $\Gamma$ and the discarded region, so the active set—and hence the surrogate domain
$\tilde{\Omega}_\ell$—grows towards $\Omega$ as $\ell$ increases; a cell active on the fine level may have no active
parent on the coarse level.

The \emph{prolongation} operator $I_{\ell-1}^\ell : \mathbb{W}_{\ell-1}^{\text{DG}} \to \mathbb{W}_\ell^{\text{DG}}$ is
the natural embedding on the background mesh, realized on each fine cell by interpolating the coarse polynomial at the
fine-level support points,
\begin{equation}
    \big(I_{\ell-1}^\ell v\big)\big|_{K_{\text{fine}}} = v\big|_{K_{\text{coarse}}}, \qquad K_{\text{fine}} \subset
    K_{\text{coarse}}.
\end{equation}
For the $p$-coarsening levels, $I_{\ell-1}^\ell$ is the analogous polynomial embedding from the lower- to the
higher-degree space on the same cell. The \emph{restriction} operator $I_\ell^{\ell-1} : \mathbb{W}_\ell^{\text{DG}}
    \to \mathbb{W}_{\ell-1}^{\text{DG}}$ is defined as the adjoint of the prolongation with respect to the $\ell^2$ inner
product on the coefficient vectors,
\begin{equation}
    I_\ell^{\ell-1} = \big(I_{\ell-1}^\ell\big)^T,
\end{equation}
which is the standard choice ensuring that the coarse-grid operator obtained by the Galerkin triple product is
symmetric whenever the fine-grid operator is. In practice these are exactly the embedding and its transpose provided by
the \texttt{MGTransfer} facilities of \texttt{deal.II}. Defining them on the background mesh sidesteps the lack of
nestedness of the surrogate spaces: the inactive degrees of freedom carried along are decoupled from the active ones
(Section~\ref{sec:implementation}) and therefore do not pollute the transferred active data.

We emphasize that the inter-level restriction $I_\ell^{\ell-1}$ introduced here is distinct from the subspace
restriction $R_i$ of Section~\ref{sec:smoother}: $R_i$ maps a function on a single level $\ell$ to the cell-local
subspace $\mathbb{V}_i$ used by the smoother, whereas $I_\ell^{\ell-1}$ transfers data between two consecutive
multigrid levels $\ell$ and $\ell-1$. We use different symbols ($R$ versus $I$) to keep the two roles clearly
separated.

We note, however, that this construction still inherits the geometric mismatch between active sets at consecutive
levels. We tried to reduce this potential efficiency loss by implementing a more sophisticated prolongation that would
extrapolate the values that are not defined on the coarser mesh, and its adjoint as a restriction. However, this
approach did not yield any significant improvements in the convergence rates. We therefore stick to the standard
prolongation and restriction operators.

The choice of the threshold parameter $\lambda$ directly influences the set of active cells, and thus the computational
domain $\tilde{\Omega}_\ell$, on each multigrid level $\ell$. A lower $\lambda$ value generally results in more cells
being classified as active, potentially leading to a surrogate domain that more closely envelops the true domain.
However, the sets of active cells across different multigrid levels, $\tilde{\Omega}_\ell$ and
$\tilde{\Omega}_{\ell-1}$, are not guaranteed to be nested. This lack of geometric nestedness can pose a challenge for
standard transfer operators. For instance, a coarse cell might be active while some of its children on the finer level
are not, or vice-versa. This mismatch can reduce the efficiency of inter-grid transfers, as information might be
prolonged to or restricted from regions that are not consistently part of the active domain across levels. While our
experiments with more sophisticated extrapolation-based transfer operators did not show significant gains, the
interplay between $\lambda$, the resulting active cell configurations, and the transfer operator effectiveness remains
an important consideration for multigrid performance in SBM.

\section{Implementation details}
\label{sec:implementation}
The numerical implementation of the DG-SBM multigrid solver is based on the open-source finite element library
\texttt{deal.II}~\cite{dealii2024}. It provides a comprehensive framework for the implementation of~finite
element methods, including mesh handling, finite element spaces, assembly of linear systems, and interfaces to various
linear algebra solvers and preconditioners. Its flexibility and extensive capabilities make it well-suited for
developing complex methods like the one presented here, which involves unfitted meshes and multigrid techniques.

The spatial discretization of the Poisson equation using the Discontinuous Galerkin method follows the standard
framework available within \texttt{deal.II}. This involves defining the DG finite element space on the background
Cartesian mesh, assembling the element-wise contributions to the global stiffness matrix and right-hand side vector,
and handling the jump and average terms across interior faces as described in Section 2.2. The library's support for
various DG formulations simplifies the implementation of the bilinear form $a_h(u_h, v_h)$.

\subsection{Background mesh and geometry handling}
When implementing SBM on a non-body-fitted mesh, we need to handle cells intersected by the true boundary $\Gamma$. We
use a level set function $\phi(\mathbf{x})$ to implicitly define the domain $\Omega = \{\mathbf{x} \mid
    \phi(\mathbf{x}) < 0\}$, with $\Gamma = \{\mathbf{x} \mid \phi(\mathbf{x}) = 0\}$. Cells of the background mesh are
classified based on their intersection with the zero level set: cells entirely inside $\Omega$ (interior), cells
entirely outside $\Omega$ (exterior), and cells intersected by $\Gamma$. Then, for the intersected cells the fraction
of the cell volume inside $\Omega$ is computed, and the cell is classified as active if this fraction is greater than a
threshold $\lambda$. This is handled using non-matching quadrature rules implemented in \texttt{deal.II}, which are
based on the techniques described in~\cite{saye2015high}.

In our approach, degrees of freedom are formally assigned to all cells of the background mesh, including those that are
classified as non-active. However, for non-active cells, the local contribution to the global system matrix is set to
the identity matrix, and the corresponding entries in the right-hand side are set to zero. This effectively decouples
the degrees of freedom in non-active cells from the system, ensuring that the solution is only computed within the
computational domain $\tilde{\Omega}$ while simplifying the global matrix assembly process by maintaining a uniform
structure of degrees of freedom (DoF) across all cells.

\subsection{Processing surrogate boundary and matrix assembly}
\label{sec:sbm_assembly}
The SBM requires computing the closest point projection from points on the surrogate boundary $\tilde{\Gamma}$ to the
true boundary $\Gamma$. For the specific test cases, where it is possible, this projection can be computed exactly.
For~
general geometries, finding the closest point involves solving a local nonlinear optimization problem for each
point on $\tilde{\Gamma}$.
The problem is to find the point on
the true boundary that minimizes the distance to a given point on the surrogate boundary, since the true boundary is
given by the zero level set of the level set function $\phi(\mathbf{x})$, this can be formulated as:
\begin{equation}
    \min_{\mathbf{x} \in \Gamma} \|\mathbf{x} - \tilde{\mathbf{x}}\|^2
    \quad \text{s.t.} \quad \phi(\mathbf{x}) = 0.
\end{equation}
We introduce a Lagrange multiplier $\mu$ to enforce the constraint $\phi(\mathbf{x}) = 0$. The Lagrangian function is:
\begin{equation}
    \mathcal{L}(\mathbf{x}, \mu) = \|\mathbf{x} - \tilde{\mathbf{x}}\|^2 + \mu \phi(\mathbf{x}).
\end{equation}
The necessary conditions for optimality are given by the stationary points of the Lagrangian, which satisfy the
following system of equations:
\begin{align}
    \nabla_{\mathbf{x}} \mathcal{L}(\mathbf{x}, \mu)           & = 2(\mathbf{x} - \tilde{\mathbf{x}}) + \mu
    \nabla\phi(\mathbf{x}) = 0,
    \\
    \frac{\partial \mathcal{L}}{\partial \mu}(\mathbf{x}, \mu) & = \phi(\mathbf{x}) = 0.
\end{align}
Solving this system yields the closest point $\mathbf{x}$ on the true boundary $\Gamma$ to the point
$\tilde{\mathbf{x}}$ on the surrogate boundary $\tilde{\Gamma}$. The nonlinear system is solved
iteratively using a Newton-Raphson method. As the Newton-Raphson method requires the evaluation of the Jacobian matrix,
we need to compute the second derivatives of the level set function $\phi(\mathbf{x})$. To ensure sufficient smoothness
for derivative calculations in the closest point search, the level set function is represented using finite elements of
degree $2$ for $p=1$
and degree $p$ for $p>1$.

We emphasize that this closest-point search is a one-time preprocessing step: the projections and the associated shift
vectors $\mathbf{d}$ are computed once, during the setup phase, for each quadrature point on the surrogate boundary
(and, in the multigrid context, once per level), and are then stored and reused; no Newton solve is repeated during the
iterative solution of the linear system. The individual solves are moreover completely independent of one another and
therefore trivially parallelizable. Since the quadrature points involved are confined to the surrogate boundary, their
number scales as $\mathcal{O}(h^{-(d-1)})$, compared to $\mathcal{O}(h^{-d})$ cells in the volume, and each Newton
iteration converges in a few steps for the smooth level sets considered here. In our experiments, the cost of this
geometric preprocessing was negligible compared to matrix assembly.

The search for the closest point is performed only in the interior of the cells adjacent to the surrogate boundary.
While this does not guarantee that the closest point is found inside the cell, it is expected that even if the closest
point is outside the cell, the extrapolation will still yield a good approximation of the boundary condition. We will
verify this assumption in the numerical results.

Finally, the extrapolation of the function values from the true boundary to the surrogate boundary required for the
matrix assembly process is accomplished by evaluating the values of the basis functions at the points on the surrogate
boundary. The assembly of the cell contributions and interior faces to the system matrix and right-hand side vector is
performed using the standard finite element assembly process provided by \texttt{deal.II}.

\subsection{Multigrid structures}
\label{sec:multigrid_structures}
For the multigrid hierarchy, we leverage \texttt{deal.II}'s built-in capabilities for handling nested meshes and
defining finite element spaces on each level. The standard projection operators provided by \texttt{deal.II} are used
for the prolongation ($I_{\ell-1}^\ell$) and restriction ($I_\ell^{\ell-1}$) operators, transferring data between
coarser and finer grid levels. These operators act on the entire background mesh, transferring the solution for all
degrees of freedom, irrespective of whether they correspond to active or non-active cells. This involves transferring
values in regions outside the computational domain. Furthermore, since the smoother and residual evaluations are
restricted to the active cells, the values in the inactive regions do not propagate into the solution within the domain
of interest, nor do they influence the convergence of the method.

The assembly of the system matrix $A_\ell$ on each level $\ell$ of the multigrid hierarchy follows a procedure
analogous to that on the finest level. The level set function defining the domain geometry is interpolated onto the
mesh $\mathcal{T}_\ell$. Based on this interpolated level set and the chosen threshold $\lambda$, active cells for
level $\ell$ are identified. The DG-SBM bilinear form is then used to assemble the local contributions to $A_\ell$ only
for these active cells. For cells deemed non-active on level $\ell$, their corresponding entries in $A_\ell$ are set to
effectively decouple them from the system by contributing an identity block to the matrix and zero to the right-hand
side. While assembling matrices on all levels is computationally expensive, assembling at least the finest level matrix
is a requirement for algebraic multigrid (AMG) methods. In contrast, geometric multigrid allows for a completely
matrix-free implementation, where level operators are evaluated on-the-fly. Recent work~\cite{wichrowski2025MFSBM} has
shown that matrix-free evaluation of SBM operators is highly efficient, suggesting that the current matrix-based
approach can be readily extended to a more performant matrix-free setting.

The smoother component of the multigrid V-cycle is implemented as a cell-wise Symmetric Successive Over-Relaxation
(SSOR) method. The SSOR iteration proceeds by visiting each cell sequentially and solving a small local linear system
corresponding to the block of degrees of freedom within that cell, using the most recently updated values from
previously processed cells. This choice is motivated by its relative ease of tuning compared to other smoothers. For
instance, a Block Jacobi smoother typically requires a carefully chosen relaxation parameter. In contrast, SSOR-based
smoothers often perform well with a fixed relaxation parameter, frequently $\omega=1$, reducing the need for
problem-specific tuning. Unless stated otherwise, we use $\omega=1$.

The current implementation focuses on testing the fundamental concepts and effectiveness of the multigrid
preconditioner for the DG-SBM method. Therefore, it is based on serial computations using sparse matrices. While
\texttt{deal.II} supports parallelization, the initial goal was to validate the multigrid algorithm's convergence
properties and identify potential challenges in the DG-SBM context before scaling up to larger, parallel problems.

\section{Numerical results}
\label{sec:numerical_results}
In this section, we present numerical results to evaluate the performance of the proposed DG-SBM multigrid
preconditioner for the Poisson equation. We investigate its effectiveness in terms of convergence rates, iteration~
counts under mesh refinement ($h$-refinement) and polynomial degree increase ($p$-refinement). To better illustrate the
multigrid preconditioner's performance characteristics and enable meaningful comparison of iteration counts, we solve
the system to high accuracy with a tolerance of $10^{-12}$ for the relative residual reduction. This stringent
tolerance results in higher iteration counts, making the differences in solver performance more apparent and
facilitating clearer assessment of the multigrid effectiveness. While a tolerance of $10^{-12}$ is stringent for many
engineering applications, we use it here to rigorously test the robustness of the multigrid preconditioner. We have
verified that using a more relaxed tolerance, such as $10^{-10}$, results in the same trends, but since the iteration
count is lower, it is harder to differentiate between different configurations. In particular, for $p=3$, the method
struggles with convergence, and the higher tolerance helps in identifying these challenges. Solver failure is defined as
exceeding 100 GMRES iterations without reaching this tolerance.

We use the relative residual reduction as the stopping criterion throughout. It is well known that the residual can
underestimate the algebraic error by a factor of the condition number, $\|u_h - u_h^k\| \le \kappa(A_h) \,
    \|r_k\|/\|A_h\|$, so that a small residual need not imply a small error for severely ill-conditioned systems. For the
present problem, however, the conditioning is benign in the sense that it is the standard one of an $H^1$-elliptic
discretization: $\kappa(A_h) = \mathcal{O}({(p+1)}^4 h^{-2})$, with the SBM extrapolation not increasing this order.
This $\mathcal{O}(h^{-2})$ growth is unfavorable for iterative solvers but is in no way pathological---it is the same
conditioning encountered in every standard finite element discretization of the Poisson problem, for which
residual-based stopping is routine. Combined with the stringent tolerance $10^{-12}$, the resulting algebraic error
stays far below the discretization error on all meshes considered. This is confirmed directly by the convergence study
in Section~\ref{sec:convergence}: because we report the error against the exact (manufactured) solution, the optimal
$\mathcal{O}(h^{p+1})$ rates observed in Figure~\ref{fig:convergence_results_l2} would not be attainable if the
algebraic error left by the residual-based criterion were polluting the solution. We therefore retain the relative
residual as a reliable and convenient stopping criterion for this class of problems. The initial background mesh is a
square (or cube in 3D) domain covering the range $[-1.01, 1.01]$ in each dimension, subdivided into $4$ cells in each
coordinate direction. We use $\sigma_\Gamma = 5$ and $\sigma_F = 1$ for the penalty parameters in the SBM formulation.
The interior penalty $\sigma_F=1$ is chosen as a minimal value to ensure coercivity while minimizing its impact on
solver convergence, whereas the boundary penalty $\sigma_\Gamma=5$ is selected to be sufficiently large to ensure
stability of the Nitsche coupling on the surrogate boundary.

\begin{figure}[h!]
    \centering
    \begin{subfigure}{0.49\textwidth}
        \includegraphics[width=\textwidth]{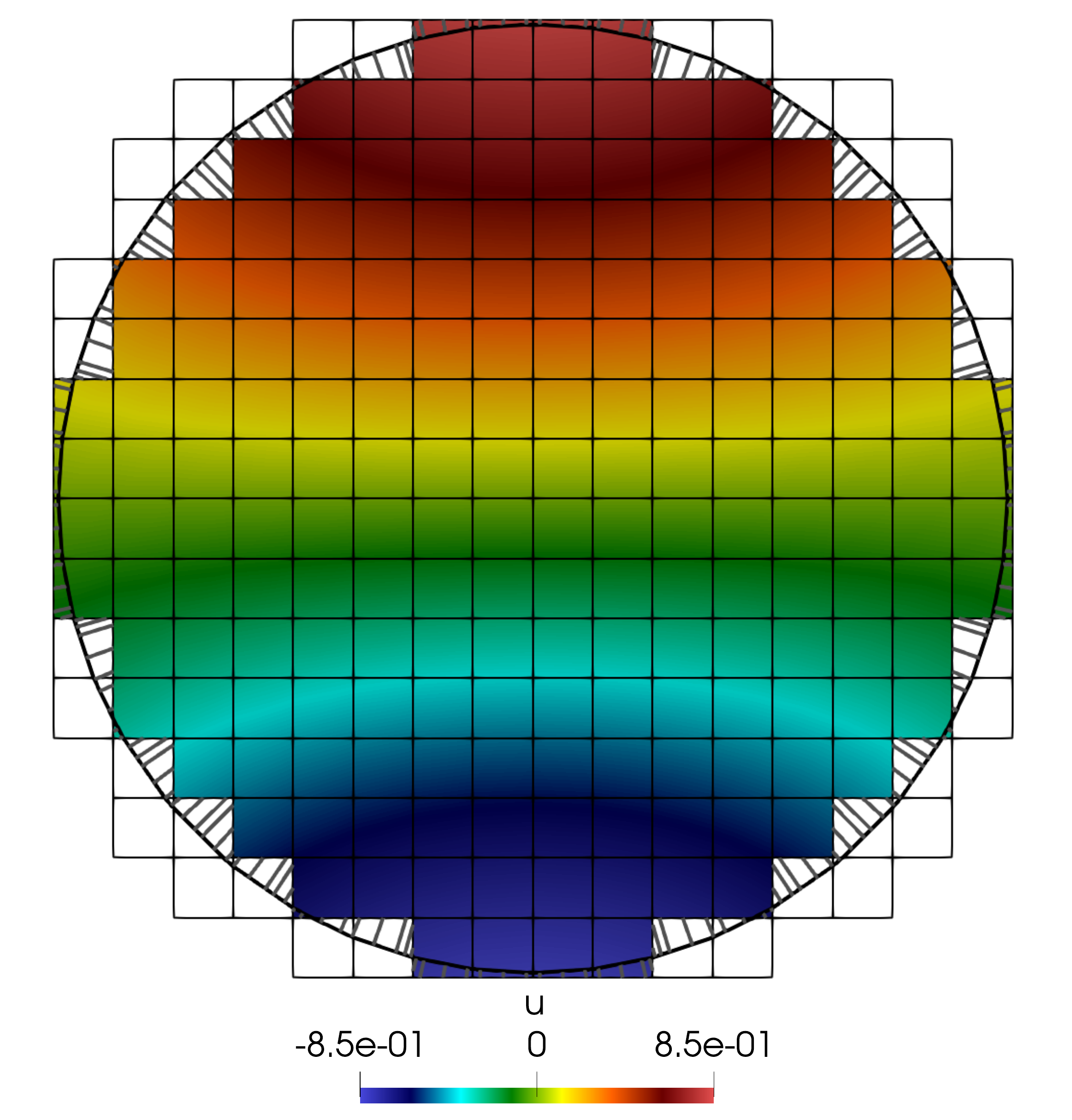}
    \end{subfigure}
    \hfill
    \begin{subfigure}{0.49\textwidth}
        \centering
        \begin{tikzpicture}
            \begin{axis}[
                    xlabel={Refinement Level},
                    ylabel={Normalized Shift Magnitude},
                    xmin=1, xmax=7,
                    xtick={1,2,3,4,5,6,7},
                    ymin=-1.5,
                    ymax=1.5,
                    ytick={-1.41421356, -1, 0,1, 1.41421356},
                    yticklabels={$-\sqrt{2}$,$-1$ ,$0$,  $1$, $\sqrt{2}$},
                    width=\textwidth, 
                    legend pos=south west, 
                    legend style={font=\scriptsize, cells={anchor=west}, legend columns=2},
                    height=7cm, 
                    grid=major,
                ]
                \addplot+[mark=*,color=red, mark options={solid}, thick,forget plot,  dashed ] coordinates {
                        (1, 0.0000) (2, 0.0000) (3, 0.0000) (4, 0.0000) (5, 0.0000) (6, 0.0000) (7, 0.0000)
                    };
                \addlegendentry{$\lambda=1.0$}
                \addplot+[mark=*,color=red, mark options={solid}, thick ] coordinates {
                        (1, 1.0825) (2, 1.1813) (3, 1.0437) (4, 1.2373) (5, 1.3157) (6, 1.3304) (7, 1.3510)
                    };

                \addplot+[mark=square*,color=darkgreen, mark options={solid}, thick, forget plot,  dashed ] coordinates
                    {
                        (1, -0.2334) (2, -0.3089) (3, -0.3904) (4, -0.3152) (5, -0.4135) (6, -0.4265) (7, -0.4302)
                    };
                \addlegendentry{$\lambda=0.75$}
                \addplot+[mark=square*,color=darkgreen, mark options={solid}, thick] coordinates {
                        (1, 0.7755) (2, 0.8005) (3, 0.7998) (4, 0.8367) (5, 0.8399) (6, 0.8450) (7, 0.8609)
                    };
                \addplot+[mark=triangle*,color=blue, mark options={solid}, thick, forget plot,	dashed ] coordinates {
                        (1, -0.4811) (2, -0.6414) (3, -0.5174) (4, -0.6021) (5, -0.6216) (6, -0.6392) (7, -0.6436)
                    };
                \addlegendentry{$\lambda=0.5$}
                \addplot+[mark=triangle*,color=blue, mark options={solid}, thick] coordinates {
                        (1, 0.3159) (2, 0.2774) (3, 0.5826) (4, 0.5807) (5, 0.6458) (6, 0.6540) (7, 0.6299)
                    };

                \addplot+[mark=diamond*,color=orange, mark options={solid}, thick, forget plot,  dashed ] coordinates {
                        (1, -0.4811) (2, -0.6414) (3, -0.7643) (4, -0.8276) (5, -0.8071) (6, -0.8421) (7, -0.8377)
                    };
                \addlegendentry{$\lambda=0.25$}
                \addplot+[mark=diamond*,color=orange, mark options={solid}, thick] coordinates {
                        (1, 0.3159) (2, 0.2774) (3, 0.2994) (4, 0.3977) (5, 0.3813) (6, 0.4077) (7, 0.4190)
                    };

            \end{axis}
        \end{tikzpicture}
    \end{subfigure}
    \caption{Left: Exemplary numerical solution for the Poisson equation on a unit disc domain, visualized on the
        surrogate domain ($\lambda=0.75$), obtained with the DG-SBM multigrid method. The true boundary $\Gamma$ is
        shown in black lines; connections between quadrature  points ($p=2$) on the surrogate boundary and their
        projections on the true boundary are illustrated with grey lines. Note that with $\lambda=0.75$,
        some parts of the true boundary lie within the surrogate domain
        Right: Illustrative distribution of minimum and maximum shift magnitudes (normalized by cell size $h$) on
        the surrogate boundary
        for the unit disk problem with $p=3$. Here the \emph{shift magnitude} is the length $|\mathbf{d}|$ of the shift
        vector $\mathbf{d} = \mathbf{x} - \tilde{\mathbf{x}}$ connecting a point $\tilde{\mathbf{x}}$ on the surrogate
        boundary to its closest-point projection $\mathbf{x}$ on the true boundary; the minimum and maximum are taken
        over all such points on the surrogate boundary at each refinement level (a signed convention is used, with
        negative values denoting shifts opposing the outward surrogate normal). The values $\pm\sqrt{2}$ represent the
        theoretical maximum possible normalized shift magnitude in 2D for a square cell.}
    \label{fig:shift_distribution}
\end{figure}

The majority of our tests are conducted on a unit sphere (a unit disk in 2D, depicted on the left panel of
Fig.~\ref{fig:shift_distribution}). While geometrically simple, the unit sphere serves as an insightful benchmark. Any
sufficiently smooth ($C^1$) complex boundary, when viewed at a fine enough mesh resolution, locally resembles a flat
plane. The unit sphere, due to its uniform curvature, presents a comprehensive range of intersection angles between the
true boundary and the background mesh cells. Furthermore, the boundary intersects cells at various locations relative
to cell centers and faces, leading to a diverse distribution of shift vector magnitudes and directions. This includes
scenarios where the shift vector points from the surrogate boundary towards the interior of the true domain (which we
denote as negative shifts if they oppose the outward normal of the surrogate boundary). The right panel of
Figure~\ref{fig:shift_distribution} depicts the minimum and maximum shift magnitudes observed across the surrogate
boundary for a typical discretization. The shift magnitudes are normalized by the cell size $h$; for a unit cell, the
theoretically largest possible shift magnitude in 2D is $\sqrt{2}h$, occurring if the surrogate boundary point is at a
cell corner and the true boundary passes through the diagonally opposite corner.

The primary metric for evaluating the multigrid preconditioner is the number of GMRES iterations required to reduce the
initial residual by a factor of $10^{-12}$. If the threshold is not met within 100 iterations, the solver is considered
to have failed. All experiments are conducted using the implementation described in the previous section. As our solver
depends on the value of the threshold parameter $\lambda$, we will first explore the solver performance on 2D problems
as they are less computationally intensive and then show some results for 3D problems with tuned $\lambda$.

\subsection{Convergence rates}
\label{sec:convergence}
We first verify convergence rates using a manufactured solution approach. For the Poisson equation on the unit sphere
$\Omega = \{\mathbf{x} \in \mathbb{R}^d \; :\; \|\mathbf{x}\| < 1\}$, we choose the manufactured solution
\begin{equation}
    u(\mathbf{x}) = 2 \cos(x_1)  \sin(x_2),
\end{equation}
from which we derive the source term $f = -\Delta u$ and Dirichlet boundary condition $g = u|_{\Gamma}$.

\begin{figure}[h!]
    \centering

    \begin{subfigure}[b]{0.49\textwidth}
        \centering
        \begin{tikzpicture}
            \begin{axis}[
                xlabel={Mesh size, $H=2.02$},
                ylabel={$L_2$ Error }, 
                xmin=1, xmax=5,
                xtick={1,2,3,4,5},
                xticklabels={${H}/{2^1}$, ${H}/{2^2}$, ${H}/{2^3}$, ${H}/{2^4}$, ${H}/{2^5}$},
                ymin=6e-11, 
                ymax=6e-2, 
                ymode=log,
                legend pos=south west, 
                legend style={font=\scriptsize, cells={anchor=west}},
                width=\textwidth,
                height=7cm,
                grid=major,
                ]

                \addplot [name path=upper_p1, mark=*,red,  mark options={solid}, thick ] coordinates {
                        (1,1.8369e-02) (2,7.3833e-03) (3,1.9051e-03) (4,6.4445e-04) (5,1.5367e-04)
                    };
                \addplot [name path=lower_p1, mark=*,red,  mark options={solid}, thick, forget plot] coordinates {
                        (1,5.1931e-03) (2,2.5188e-03) (3,4.2610e-04) (4,8.3951e-05) (5,1.7045e-05)
                    };
                \addlegendentry{$p=1$}
                \addplot [fill=red, fill opacity=0.3,forget plot] fill between [of=upper_p1 and lower_p1];

                \addplot [name path=upper_p2, mark=square*,darkgreen, mark options={solid}, thick, forget plot]
                coordinates {
                        (1,2.0915e-03) (2,2.2065e-004) (3,2.0348e-05) (4,2.6120e-06) (5,2.9301e-07)
                    };
                \addplot [name path=lower_p2, mark=square*,darkgreen, mark options={solid}, thick ]
                coordinates {
                        (1,1.2167e-03) (2,1.5255e-04) (3,1.2272e-05) (4,9.6671e-07) (5,8.8114e-08)
                    };
                \addlegendentry{$p=2$}
                \addplot [fill=darkgreen, fill opacity=0.3, forget plot] fill between [of=upper_p2 and lower_p2];

                \addplot [name path=upper_p3, mark=triangle*,blue,	mark options={solid}, thick, forget plot]
                coordinates {
                        (1,1.3594e-04) (2,1.2242e-05) (3,6.6510e-07) (4,6.7293e-08) (5,4.1223e-09)
                    };
                \addplot [name path=lower_p3, mark=triangle*,blue,	mark options={solid}, thick] coordinates {
                        (1,4.0585e-05) (2,4.7205e-06) (3,1.9789e-07) (4,9.2215e-09) (5,4.3084e-10)
                    };
                \addlegendentry{$p=3$}
                \addplot [fill=blue, fill opacity=0.3, forget plot] fill between [of=upper_p3 and lower_p3];


                \addplot [dashed, red, thick, forget plot] coordinates {
                        (1,1.7e-02)
                        (5,1.7e-02 / 256)
                    };
                \addplot [dashed, darkgreen, thick, forget plot] coordinates {
                        (1,2.5e-03)
                        (5,2.5e-03 / 4096)
                    };
                \addplot [dashed, blue, thick, forget plot] coordinates {
                        (1,1.e-04)
                        (5,1.e-04 / 65536)
                    };
            \end{axis}
        \end{tikzpicture}
    \end{subfigure}
    \caption{Influence of threshold $\lambda$ on $L_2$ error versus mesh size $h = \frac{H}{2^n}$ for the DG-SBM with $p=1, 2, 3$ on the unit sphere domain
        in 2D. Shaded regions show the range of $L_2$ errors for $\lambda \in \{0.25, 0.5, 0.75, 1.0\}$. Dashed lines
        indicate
        theoretical optimal convergence rates.}
    \label{fig:convergence_results_l2}
\end{figure}
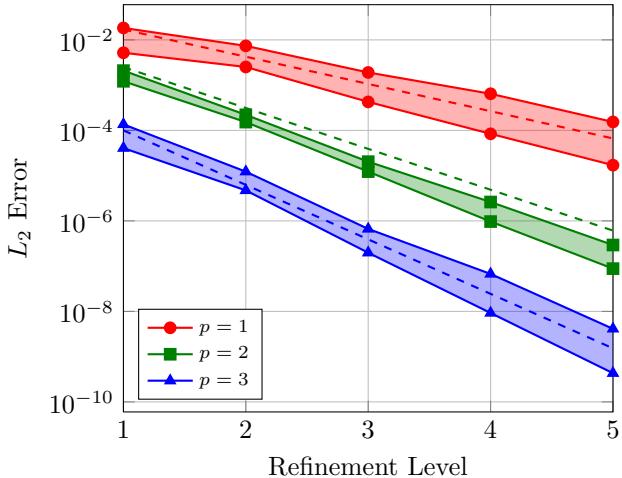

For a DG method using piecewise polynomials of degree $p$, we expect to observe optimal convergence rates of order
$\mathcal{O}(h^{p+1})$ in the $L_2$ norm, provided the solution is sufficiently smooth. Figure
\ref{fig:convergence_results_l2} shows the $L_2$ error versus the mesh size $h = 2^{-n}H$ (where $n$ is the refinement
level and $H=2.02$ is the initial mesh size) for different polynomial degrees $p$ in 2D. The shaded regions in the 2D
plot (Figure \ref{fig:convergence_results_l2}) represent the range of $L_2$ errors obtained for different values of the
threshold parameter $\lambda \in \{1.0, 0.75, 0.5, 0.25\}$, which determines the active cells. It was demonstrated
in~\cite{yang2024optimal} that the SBM may not yield stable convergence rates for $\lambda=0$. Within this range, lower
values of $\lambda$ typically lead to lower $L_2$ errors.

The results demonstrate that the DG-SBM method achieves the expected optimal convergence rates of order
$\mathcal{O}(h^{p+1})$ in the $L_2$ norm for the 2D case across these $\lambda$ values. The dashed lines, representing
theoretical convergence rates, align well with the computed errors for polynomial degrees $p=1$ through $p=3$.

\subsection{Multigrid for linear elements in 2D}
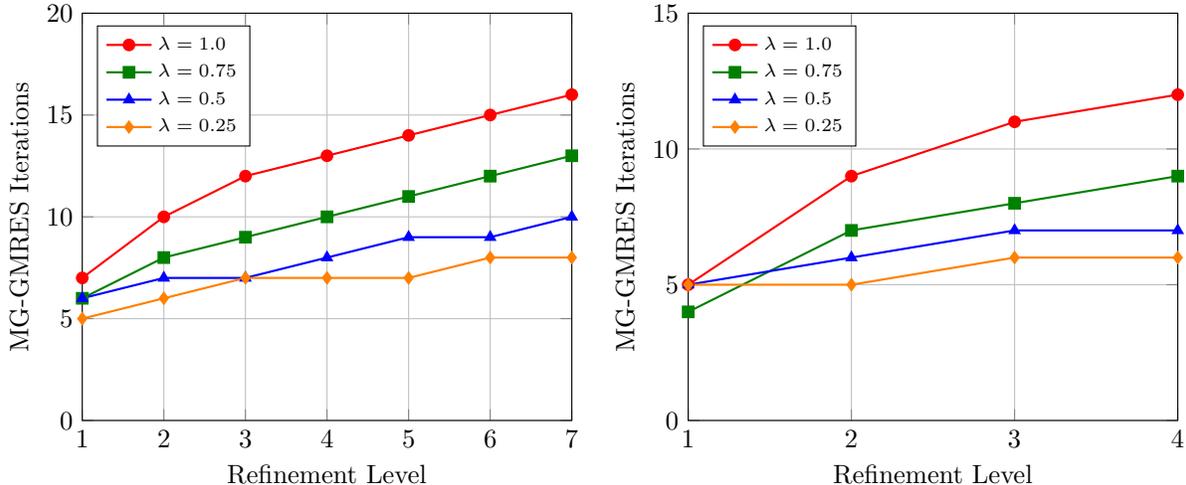
\begin{figure}[h!]
    \centering
    \begin{subfigure}[b]{0.49\textwidth}
        \centering

        \begin{tikzpicture}
            \begin{axis}[
                    xlabel={Refinement Level },
                    ylabel={MG-GMRES Iterations},
                    xmin=1, xmax=7, 
                    xtick={1,2,3,4,5,6,7},
                    ymin=0, 
                    ymax=20, 
                    legend pos=north west,
                    legend style={font=\scriptsize, cells={anchor=west}},
                    width=\textwidth, 
                    height=7cm, 
                    grid=major,
                ]
                \addplot+[mark=*,red,  mark options={solid}, thick] coordinates {(1,7) (2,10) (3,12) (4,13) (5,14)
                        (6,15) (7,16) (8,17)};
                \addlegendentry{$\lambda=1.0$}
                \addplot+[mark=square*,darkgreen,	mark options={solid}, thick] coordinates {(1,6) (2,8) (3,9)
                        (4,10) (5,11) (6,12) (7,13)};
                \addlegendentry{$\lambda=0.75$}
                \addplot+[mark=triangle*,blue,	mark options={solid}, thick] coordinates {(1,6) (2,7) (3,7)
                        (4,8) (5,9) (6,9) (7,10) (8,11)};
                \addlegendentry{$\lambda=0.5$}
                \addplot+[mark=diamond*,orange,  mark options={solid}, thick] coordinates {(1,5) (2,6) (3,7) (4,7)
                        (5,7) (6,8) (7,8) (8,8)};
                \addlegendentry{$\lambda=0.25$}

            \end{axis}
        \end{tikzpicture}

    \end{subfigure}
    \begin{subfigure}[b]{0.49\textwidth}
        \begin{tikzpicture}
            \begin{axis}[
                    xlabel={Refinement Level }, ylabel={MG-GMRES Iterations}, xmin=1, xmax=4, xtick={1,2,3,4}, ymin=0, ymax=15,
                    legend pos=north west, legend style={font=\scriptsize, cells={anchor=west}}, width=\textwidth, height=7cm, grid=major,
                ]
                \addplot+[mark=*,red,  mark options={solid}, thick] coordinates {
                        (1,5) (2,9) (3,11) (4,12)
                    };
                \addlegendentry{$\lambda=1.0$}

                \addplot+[mark=square*,darkgreen,  mark options={solid}, thick] coordinates {
                        (1,4) (2,7) (3,8) (4,9)
                    };
                \addlegendentry{$\lambda=0.75$}

                \addplot+[mark=triangle*,blue,	mark options={solid}, thick] coordinates {
                        (1,5) (2,6) (3,7) (4,7)
                    };
                \addlegendentry{$\lambda=0.5$}

                \addplot+[mark=diamond*,orange,  mark options={solid}, thick] coordinates {
                        (1,5) (2,5) (3,6) (4,6)
                    };
                \addlegendentry{$\lambda=0.25$}

            \end{axis}
        \end{tikzpicture}
    \end{subfigure}
    \caption{Convergence of GMRES: iteration counts with $h$-refinement for the SBM with different thresholds $\lambda$
        and
        $p=1$	     Left: 2D. Right: 3D. }
    \label{fig:convergence_mg_h}
\end{figure}

An efficient multigrid preconditioner should ideally exhibit $h$-independence, where the number of iterations for
convergence remains relatively constant as the mesh is refined. We investigate this for the Poisson problem on the unit
sphere. Figure~\ref{fig:convergence_mg_h} shows the number of MG-GMRES iterations versus mesh refinement for linear
elements ($p=1$). The left panel displays results for the 2D case, while the right panel shows the corresponding
results for the 3D case, both for various threshold values $\lambda$.

For $p=1$ in 2D (left panel), the iteration counts show a mild increase with mesh refinement for all tested $\lambda$
values. For $\lambda=1.0$, iterations increase from 7 to 16 as the mesh is refined from level 1 to 7. Smaller $\lambda$
values generally lead to lower and more stable iteration counts; for instance, with $\lambda=0.25$, iterations range
from 5 to 8. In 3D (right panel), a similar trend is observed for $p=1$. Iteration counts increase mildly with
refinement. For $\lambda=1.0$, iterations go from 5 to 12 across levels 1 to 4. Again, smaller $\lambda$ values like
$\lambda=0.25$ yield more stable and lower iteration counts, ranging from 5 to 6. These results suggest that for linear
elements, the multigrid preconditioner provides good scalability with mesh refinement in both 2D and 3D, with lower
$\lambda$ values (which include more intersected cells) generally leading to better solver performance.

We observe that the iteration counts for the MG-GMRES solver show low mesh dependence, especially for low values of
$\lambda$, despite the use of global transfer operators that act on the entire background mesh. This suggests that the
degrees of freedom located in the inactive regions, which are decoupled from the active ones as discussed in
Section~\ref{sec:multigrid_structures}, do not significantly hinder the convergence of the multigrid method.

\subsection{Higher order elements}
For polynomial degrees $p > 1$, we employ a $p$-multigrid strategy. This approach uses the finest geometric mesh with
the desired polynomial degree $p$ and creates additional multigrid levels by systematically reducing the polynomial
degree down to $p=1$. Since the underlying geometric mesh remains the same across these $p$-multigrid levels, the same
set of active cells is considered on each level, ensuring geometric consistency throughout the hierarchy. This
$p$-coarsening ensures that the finite element spaces are properly nested across multigrid levels, which allows us to
test the effectiveness of the smoother.

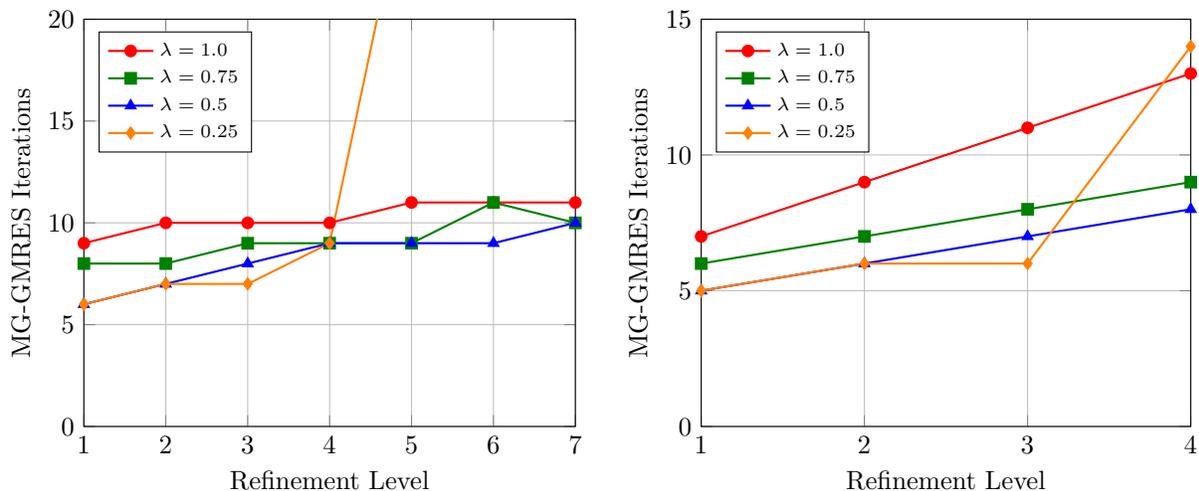
\begin{figure}[h!]
    \centering

    \begin{subfigure}[b]{0.49\textwidth}
        \begin{tikzpicture}
            \begin{axis}[
                    xlabel={Refinement Level},
                    ylabel={MG-GMRES Iterations},
                    xmin=1, xmax=7,
                    ymin=0, ymax=20,
                    xtick={1,2,3,4,5,6,7},
                    legend pos=north west,
                    legend style={font=\scriptsize, cells={anchor=west}},
                    height=7cm,
                    grid=major,
                ]
                \addplot+[mark=*, color=red, mark options={solid}, thick] coordinates {
                        (1, 9)
                        (2, 10)
                        (3, 10)
                        (4, 10)
                        (5, 11)
                        (6, 11)
                        (7, 11)
                    };
                \addlegendentry{$\lambda=1.0$}

                \addplot+[mark=square*, color=darkgreen, mark options={solid}, thick] coordinates {
                        (1, 8)
                        (2, 8)
                        (3, 9)
                        (4, 9)
                        (5, 9)
                        (6, 11)
                        (7, 10)
                    };
                \addlegendentry{$\lambda=0.75$}

                \addplot+[mark=triangle*, color=blue, mark options={solid}, thick] coordinates {
                        (1, 6)
                        (2, 7)
                        (3, 8)
                        (4, 9)
                        (5, 9)
                        (6, 9)
                        (7, 10)
                    };
                \addlegendentry{$\lambda=0.5$}

                \addplot+[mark=diamond*, color=orange, mark options={solid}, thick] coordinates {
                        (1, 6)
                        (2, 7)
                        (3, 7)
                        (4, 9)
                        (5, 28)
                        (6, 100)
                    };
                \addlegendentry{$\lambda=0.25$}
            \end{axis}
        \end{tikzpicture}

    \end{subfigure}
    \begin{subfigure}[b]{0.49\textwidth}
        \begin{tikzpicture}
            \begin{axis}[
                    xlabel={Refinement Level }, ylabel={MG-GMRES Iterations}, xmin=1, xmax=4, xtick={1,2,3,4,5}, ymin=0, ymax=15,
                    legend pos=north west, legend style={font=\scriptsize, cells={anchor=west}}, width=\textwidth, height=7cm, grid=major,
                ]
                \addplot+[mark=*,red,  mark options={solid}, thick] coordinates {
                        (1,7) (2,9) (3,11) (4,13)
                    };
                \addlegendentry{$\lambda=1.0$}

                \addplot+[mark=square*,darkgreen,  mark options={solid}, thick] coordinates {
                        (1,6) (2,7) (3,8) (4,9)
                    };
                \addlegendentry{$\lambda=0.75$}

                \addplot+[mark=triangle*,blue,	mark options={solid}, thick] coordinates {
                        (1,5) (2,6) (3,7) (4,8)
                    };
                \addlegendentry{$\lambda=0.5$}

                \addplot+[mark=diamond*,orange,  mark options={solid}, thick] coordinates {
                        (1,5) (2,6) (3,6) (4,14)
                    };
                \addlegendentry{$\lambda=0.25$}

            \end{axis}
        \end{tikzpicture}
    \end{subfigure}
    \caption{Convergence of MG-GMRES solver: iteration counts versus mesh refinement level for the DG-SBM method on the
        unit sphere domain. Results are shown for quadratic elements ($p=2$) with different threshold values $\lambda
            \in \{0.25, 0.5, 0.75, 1.0\}$. Left: 2D case. Right: 3D case.}
    \label{fig:convergence_3d_p1_p2_sphere}
\end{figure}

Figure~\ref{fig:convergence_3d_p1_p2_sphere} presents the GMRES iteration counts for polynomial degree $p=2$ under mesh
refinement for the unit sphere problem in 2D (left panel) and 3D (right panel). In 2D, the iteration counts for $p=2$
show low mesh dependence, particularly for $\lambda=0.75$ and $\lambda=1.0$, where iterations remain around 8--11. For
$\lambda=0.5$, iterations are also stable, increasing slightly from 6 to 10. However, for $\lambda=0.25$, the solver
performance degrades significantly at finer refinement levels, with iteration counts jumping to 28 at level 5 and
failing at level 6. This highlights the sensitivity to $\lambda$ for $p=2$ even in 2D. In 3D, for $p=2$, the iteration
counts generally show a rease with refinement for $\lambda=0.75$ (6--9 iterations) and $\lambda=0.5$ (5--8 iterations).
For $\lambda=1.0$, iterations increase from 7 to 13. The case $\lambda=0.25$ shows stable low iterations (5--6) up to
level 3, but then jumps to 14 iterations at level 4. Overall, the $hp$-multigrid preconditioner performs reasonably
well for $p=2$, though the choice of $\lambda$ remains important.

\begin{figure}[h!]
    \centering

    \begin{subfigure}[b]{0.49\textwidth}
        \begin{tikzpicture}
            \begin{axis}[
                    xlabel={Refinement Level},
                    ylabel={MG-GMRES Iterations},
                    xmin=1, xmax=7,
                    ymin=0, ymax=60,
                    xtick={1,2,3,4,5,6,7},
                    legend pos=north west,
                    legend style={font=\scriptsize, cells={anchor=west}},
                    height=7cm,
                    grid=major,
                ]
                \addplot+[mark=*, color=red, mark options={solid}, thick] coordinates {
                        (1, 10)
                        (2, 11)
                        (3, 11)
                        (4, 13)
                        (5, 21)
                        (6, 28)
                        (7, 57)
                    };
                \addlegendentry{$\lambda=1.0$}

                \addplot+[mark=square*, color=darkgreen, mark options={solid}, thick] coordinates {
                        (1, 10)
                        (2, 9)
                        (3, 11)
                        (4, 11)
                        (5, 12)
                        (6, 12)
                        (7, 14)
                    };
                \addlegendentry{$\lambda=0.75$}

                \addplot+[mark=triangle*, color=blue, mark options={solid}, thick] coordinates {
                        (1, 8)
                        (2, 7)
                        (3, 10)
                        (4, 11)
                        (5, 10)
                        (6, 43)
                        (7, 78)
                    };
                \addlegendentry{$\lambda=0.5$}

                \addplot+[mark=diamond*, color=orange, mark options={solid}, thick] coordinates {
                        (1, 8)
                        (2, 7)
                        (3, 27)
                        (4, 18)
                        (5, 100)
                    };
                \addlegendentry{$\lambda=0.25$}
            \end{axis}
        \end{tikzpicture}

    \end{subfigure}
    \begin{subfigure}[b]{0.49\textwidth}
        \begin{tikzpicture}
            \begin{axis}[
                    xlabel={Refinement Level }, ylabel={MG-GMRES Iterations}, xmin=1, xmax=3, xtick={1,2,3}, ymin=0, ymax=70,
                    legend pos=north east, legend style={font=\scriptsize, cells={anchor=west}}, width=\textwidth, height=7cm, grid=major,
                ]
                \addplot+[mark=*,red,  mark options={solid}, thick] coordinates {
                        (1,12) (2,100)
                    };
                \addlegendentry{$\lambda=1.0$}

                \addplot+[mark=square*,darkgreen,  mark options={solid}, thick] coordinates {
                        (1,43) (2,11) (3,24)
                    };
                \addlegendentry{$\lambda=0.75$}

                \addplot+[mark=triangle*,blue,	mark options={solid}, thick] coordinates {
                        (1,9) (2,11) (3,51)
                    };
                \addlegendentry{$\lambda=0.5$}

                \addplot+[mark=diamond*,orange,  mark options={solid}, thick] coordinates {
                        (1,16)
                        (2,100)
                    };
                \addlegendentry{$\lambda=0.25$}

            \end{axis}
        \end{tikzpicture}
    \end{subfigure}
    \caption{Convergence of GMRES: iteration counts for the SBM with different thresholds $\lambda$ on a sequence of
        successively refined meshes, for polynomial degree $p=3$, with relaxation parameter $\omega=0.8$.  Left: 2D,
        Right: in 3D.}
    \label{fig:convergence_3d_p3_sphere}
\end{figure}
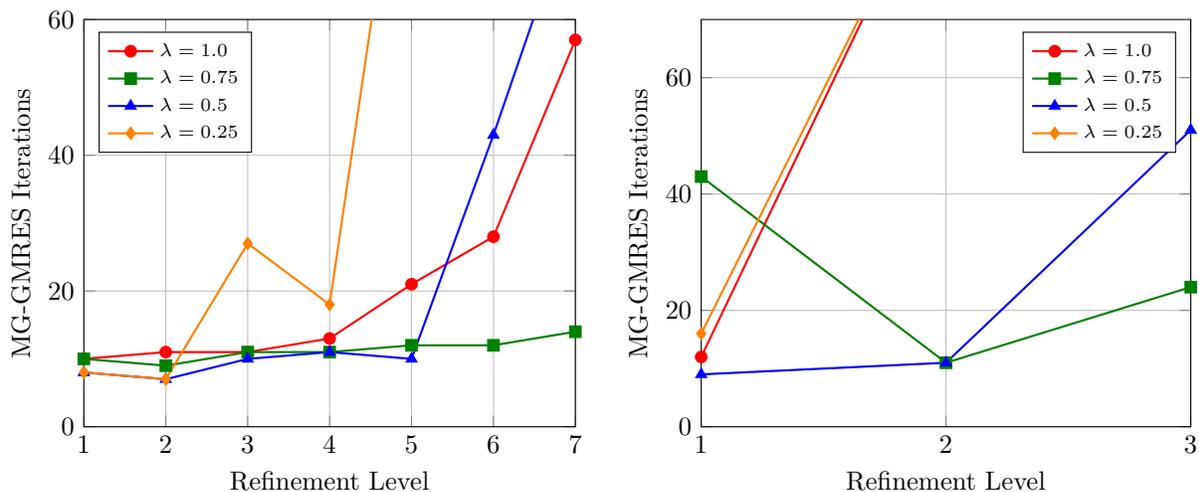

The performance for cubic elements ($p=3$) is shown in Figure~\ref{fig:convergence_3d_p3_sphere}. For these tests, the
relaxation parameter for the cell-wise SSOR smoother was set to $\omega=0.8$, as the standard $\omega=1.0$ often led to
divergence. In 2D (left panel), for $p=3$ with $\omega=0.8$, the behavior is mixed. For $\lambda=0.75$, iteration
counts are relatively stable, increasing from 10 to 14. However, for $\lambda=1.0$, iterations grow significantly with
refinement, reaching 57 at level 7. For $\lambda=0.5$, after an initial dip, iterations increase to 43 and then 78. For
$\lambda=0.25$, the solver fails at level 5. This indicates increased difficulty for $p=3$ in 2D. In 3D (right panel),
the challenges with $p=3$ are more pronounced. For $\lambda=1.0$ and $\lambda=0.25$, the solver fails (100 iterations)
already at the second refinement level. For $\lambda=0.75$, iterations are erratic (43, 11, 24). For $\lambda=0.5$,
iterations increase from 9 to 11, then jump to 51 at level 3. These results suggest that the current smoother is less
effective for $p=3$, especially in 3D, even with a tuned relaxation parameter. Due to prohibitive memory requirements
for $p \geq 3$, computations on finer grids for $p=3$ in 3D were not feasible, making it difficult to draw definitive
conclusions about the solver's scalability to finer meshes at this polynomial degree.

Numerous attempts were made to improve the convergence for $p=3$ elements, particularly in 3D. These included exploring
additive cell-wise smoothers (block Jacobi) as an alternative to the successive (block Gauss-Seidel) approach.
Furthermore, to address the challenge of selecting an optimal relaxation parameter for these additive smoothers, they
were wrapped within a Chebyshev polynomial smoother. The Chebyshev smoother can automatically determine a suitable
relaxation schedule based on estimates of the estimate of eigenvalues of the preconditioned operator. Despite these
efforts, significant improvements in iteration counts for $p=3$ cases, especially in 3D, were not consistently
observed.

\subsection{Complex geometry test case}
While the unit sphere provides a valuable test case for verifying convergence rates and basic multigrid performance, it
is also important to evaluate the method on problems with more complex geometries, where the SBM and the multigrid
preconditioner face greater challenges.

Here, we present results for the Poisson equation on a 2D domain with a non-spherical, deformed boundary. This serves
as a more challenging test for the SBM, particularly for the closest point projection algorithm which, for such
geometries, relies on the local nonlinear solver described in Section~\ref{sec:sbm_assembly}.

\begin{figure}[h!]
    \centering
    \begin{subfigure}{0.49\textwidth}
        \includegraphics[width=\textwidth]{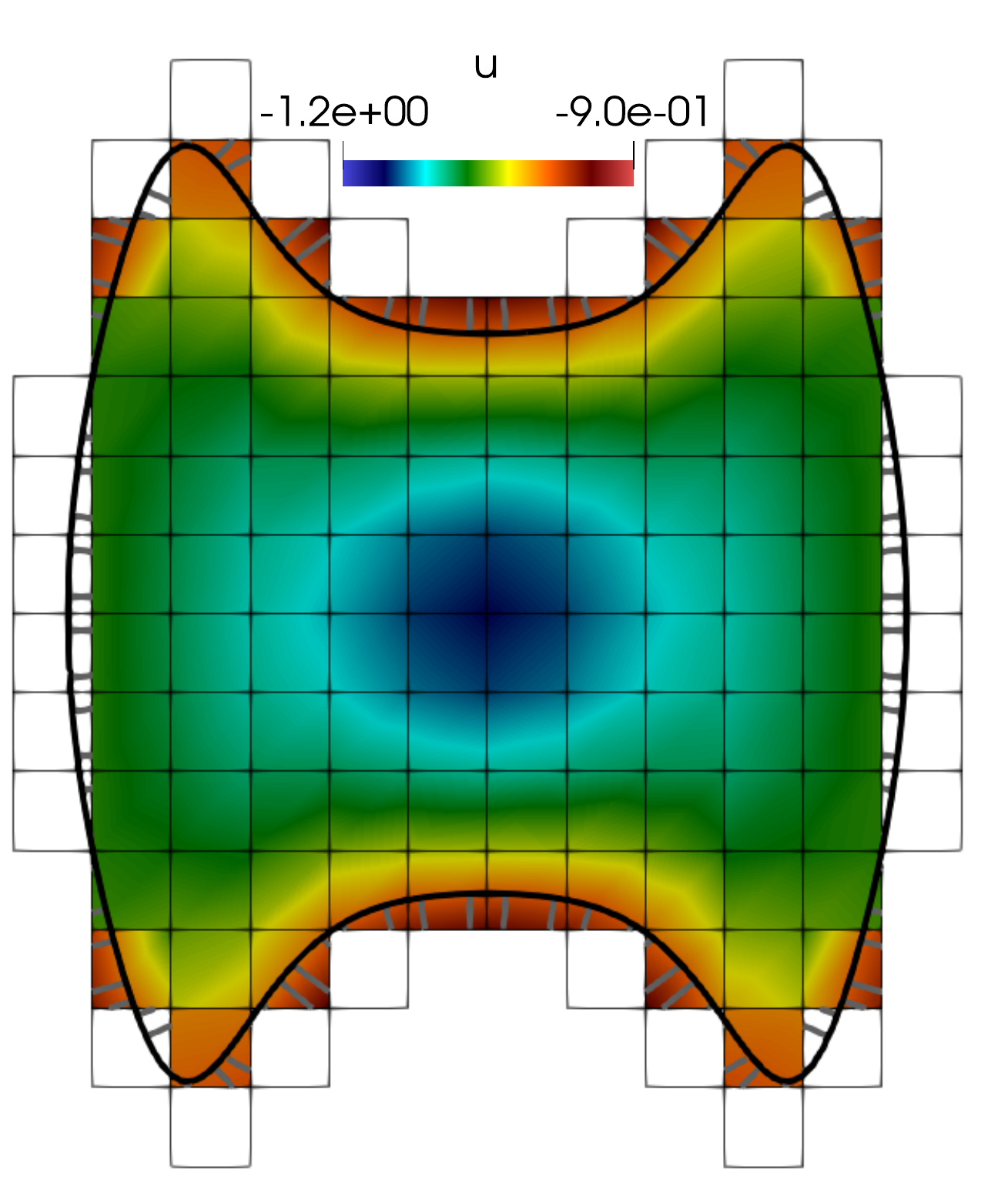}
    \end{subfigure}
    \hfill
    \begin{subfigure}{0.49\textwidth}
        \centering
        \begin{tikzpicture}
            \begin{axis}[
                xlabel={Mesh size, $H = 2.02$},
                ylabel={$L_2$ Error }, 
                xmin=1, xmax=7,
                xtick={1,2,3,4,5,6,7},
                xticklabels={${H}/{2^{1}}$, ${H}/{2^{2}}$, ${H}/{2^{3}}$, ${H}/{2^{4}}$, ${H}/{2^{5}}$, ${H}/{2^{6}}$, ${H}/{2^{7}}$},
                ymode=log, legend pos=south west, legend style={font=\scriptsize, cells={anchor=west}}, width=\textwidth, height=8cm,
                grid=major,
                ]

                \addplot [name path=upper_p1, mark=*,red,  mark options={solid}, thick ] coordinates {
                        (1,8.4530e-02) (2,4.0334e-02) (3,1.5593e-02) (4,4.5108e-03) (5,1.1889e-03) (6,3.4161e-04)
                        (7,8.5034e-05)
                    };
                \addplot [name path=lower_p1, mark=*,red,  mark options={solid}, thick, forget plot]
                coordinates {
                        (1,3.5541e-02) (2,1.6628e-02) (3,3.9802e-03) (4,1.0725e-03) (5,1.7365e-04) (6,3.3670e-05)
                        (7,6.0549e-06)
                    };
                \addlegendentry{$p=1$}
                \addplot [fill=red, fill opacity=0.3,forget plot] fill between [of=upper_p1 and lower_p1];

                \addplot [name path=upper_p2, mark=square*,darkgreen, mark options={solid}, thick,
                    forget plot]
                coordinates {
                        (1, 1.0187e-01) (2, 1.1265e-02) (3, 7.8503e-04) (4, 9.7274e-05) (5, 9.7857e-06)
                        (6, 1.2107e-06) (7, 1.4127e-07)
                    };
                \addplot [name path=lower_p2, mark=square*,darkgreen, mark options={solid}, thick ]
                coordinates {
                        (1, 6.1130e-02) (2, 1.1265e-02) (3, 5.3369e-04) (4, 4.2330e-05) (5, 3.0870e-06)
                        (6, 3.2390e-07) (7, 3.6923e-08)
                    };
                \addlegendentry{$p=2$}
                \addplot [fill=darkgreen, fill opacity=0.3, forget plot] fill between [of=upper_p2 and
                        lower_p2];

                \addplot [name path=upper_p3, mark=triangle*,blue,	mark options={solid},
                    thick, forget plot]
                coordinates {
                        (1, 2.8620e-02) (2, 1.3947e-02) (3, 3.4363e-03) (4, 2.5002e-05) (5, 1.3327e-06) (6, 1.0362e-07)
                        (7, 1.0838e-08 )
                    };
                \addplot [name path=lower_p3, mark=triangle*,blue,	mark options={solid},
                    thick] coordinates {
                        (1, 9.7938e-03) (2, 6.5175e-03) (3, 1.6143e-04) (4, 1.3241e-05) (5, 9.1573e-07)
                        (6, 7.5701e-08) (7, 2.2816e-09)
                    };
                \addlegendentry{$p=3$}
                \addplot [fill=blue, fill opacity=0.3, forget plot] fill between [of=upper_p3
                        and lower_p3];

                \addplot [dashed, red, thick, forget plot] coordinates {
                        (1,1.7e-01)
                        (7,1.7e-01 / 256/16)
                    };
                \addplot [dashed, darkgreen, thick, forget plot] coordinates {
                        (1,2.5e-02)
                        (7,2.5e-02 / 4096 / 64)
                    };
                \addplot [dashed, blue, thick, forget plot] coordinates {
                        (1,4.e-02)
                        (7,4.e-02 / 65536/256)
                    };
            \end{axis}
        \end{tikzpicture}
    \end{subfigure}
    \caption{ Left: Numerical solution for the Poisson equation on a complex deformed domain, visualized on the
        surrogate domain ($\lambda=0.5$), obtained with the DG-SBM multigrid method. The true boundary $\Gamma$ is
        shown in
        black lines; connections between quadrature points on the surrogate boundary and their projections on the true
        boundary
        are illustrated with grey lines. Note that with $\lambda=0.5$, some parts of the true boundary may lie within
        the
        surrogate domain. Right: Convergence of the DG-SBM method versus mesh size  on the same domain.}
    \label{fig:solution_complex}
\end{figure}

Our test domain $\Omega$ is defined by the zero level set of the function
\begin{equation}
    \phi(x,y) = \left( 1 - \frac{3}{4} \sin(\pi x)^2 \right) (x^2 + y^2) - \frac{1}{5},
\end{equation}
where $\Omega = \{\mathbf{x} \mid \phi(\mathbf{x}) < 0\}$. The visualization (Figure~\ref{fig:solution_complex}, left
panel) illustrates all cells of the background mesh that are at least partially inside $\Omega$. However, only those
cells where the fraction of their area inside $\Omega$ exceeds the threshold $\lambda=0.5$ are considered active. The
solution is consequently visualized on this set of active cells, which constitute the surrogate domain
$\tilde{\Omega}$. It is worth noting that even for the mesh depicted, which is not the coarsest level used in
computations, the cell diameter can be larger than the smallest radius of curvature of the domain. Consequently, once
can expect that the geometry will not be fully resolved on the coarsest multigrid levels. The DG-SBM method handles the
intricate boundary features, and the grey lines depict the shifts from the surrogate boundary to the true boundary.

To assess the convergence behavior of the DG-SBM method, particularly in the presence of a complex, implicitly defined
boundary, we employ a manufactured solution directly related to the level set function. We define the exact solution as
$u(\mathbf{x}) = \phi(\mathbf{x}) + 1$. Consequently, the source term for the Poisson equation is $f = -\Delta u
    =-\Delta \phi$. The Dirichlet boundary condition on $\Gamma = \{\mathbf{x} \mid \phi(\mathbf{x}) = 0\}$ is $g =
    u|_{\Gamma} = 1$.

This choice of a constant boundary condition is particularly insightful. Any inaccuracies in the geometric
representation of the boundary $\Gamma$ (e.g., due to approximating the level set function $\phi$ with finite elements)
will result in the constant boundary condition $g=1$ being enforced at slightly perturbed locations. These geometric
perturbations directly influence the computed solution $u_h$ and are thus reflected in the observed convergence rates.

We compute the $L_2$ error $\|u - u_h\|_{L_2(\tilde{\Omega})}$ on the surrogate domain $\tilde{\Omega}$ for a sequence
of uniformly refined meshes, and the results of this convergence study are presented in the right panel of
Figure~\ref{fig:solution_complex}. We note that for $p=3$ the MG-GMRES solver failed to converge for some $\lambda$
values, leading to a lack of data for these cases. The shaded region for $p=3$ is constructed using data from the
subset of $\lambda$ values for which the solver converged at each respective refinement level.

The plot shows that while optimal convergence rates are generally achieved on finer meshes, convergence on coarser
meshes, especially for $p=3$, can initially appear stagnant. This occurs because the geometry is under-resolved by the
coarser discretizations. Once the cell size is small enough to resolve the geometry, the optimal rate is achieved.
However, when the cell size is comparable to or larger than the local radius of curvature, the discrete problem
effectively addresses a smoothed or simplified version of the true boundary. This can lead to an initial, but
unreliable, rapid error reduction until the mesh becomes fine enough to capture the actual geometric details.

\begin{figure}[h!]
    \centering
    \begin{subfigure}{0.49\textwidth}

        \begin{tikzpicture}
            \begin{axis}[
                    xlabel={Refinement Level },
                    ylabel={MG-GMRES Iterations},
                    xmin=1, xmax=7,
                    xtick={1,2,3,4,5,6,7},
                    ymin=0,
                    ymax=45, 
                    legend pos=north west,
                    legend style={font=\scriptsize, cells={anchor=west}},
                    width=\textwidth,
                    height=7cm,
                    grid=major,
                ]
                \addplot+[mark=*,red,  mark options={solid}, thick] coordinates {
                        (1,3) (2,10) (3,13) (4,15) (5,16) (6,17) (7,17)
                    };
                \addlegendentry{$\lambda=1.0$}

                \addplot+[mark=square*,darkgreen,	mark options={solid}, thick] coordinates {
                        (1,4) (2,10) (3,11) (4,12) (5,12) (6,13) (7,14)
                    };
                \addlegendentry{$\lambda=0.75$}

                \addplot+[mark=triangle*,blue,	mark options={solid}, thick] coordinates {
                        (1,5) (2,7) (3,8) (4,8) (5,9) (6,9) (7,10)
                    };
                \addlegendentry{$\lambda=0.5$}

                \addplot+[mark=diamond*,orange,  mark options={solid}, thick] coordinates {
                        (1,5) (2,6) (3,7) (4,7) (5,8) (6,8) (7,8)
                    };
                \addlegendentry{$\lambda=0.25$}

                \addplot+[mark=*,red,  mark options={solid}, thick, dashed] coordinates {
                        (1,11) (2,17) (3,20) (4,24) (5,27) (6,29) (7,30)
                    };

                \addplot+[mark=square*,darkgreen,	mark options={solid}, thick, dashed] coordinates {
                        (1,9) (2,14) (3,18) (4,20) (5,22) (6,23) (7,24)
                    };

                \addplot+[mark=triangle*,blue,	mark options={solid}, thick, dashed] coordinates {
                        (1,9) (2,11) (3,13) (4,14) (5,15) (6,15) (7,16)
                    };

                \addplot+[mark=diamond*,orange,  mark options={solid}, thick, dashed] coordinates {
                        (1,7) (2,9) (3,13) (4,12) (5,12) (6,13) (7,13)
                    };

                \addplot+[mark=*,red,  mark options={solid}, thick, dotted] coordinates {
                        (1,14) (2,16) (3,23) (4,29) (5,34) (6,46)(7, 100)
                    };

                \addplot+[mark=square*,darkgreen,	mark options={solid}, thick, dotted] coordinates {
                        (1,15) (2,15) (3,20) (4,24) (5,47) (6,29) (7,38)
                    };

                \addplot+[mark=triangle*,blue,	mark options={solid}, thick, dotted] coordinates {
                        (1,9) (2,7) (3,17) (4,26) (5,23) (6,23)(7,100)
                    };

                \addplot+[mark=diamond*,orange,  mark options={solid}, thick, dotted] coordinates {
                        (1,13) (2,10) (3,19)(4,100)
                    };

            \end{axis}
        \end{tikzpicture}
    \end{subfigure}

    \caption{Convergence of GMRES: iteration counts with multigrid preconditioner for the DG-SBM method on a complex
        deformed
        domain in 2D. Solid lines: $p=1$, dashed: $p=2$, dotted: $p=3$. The results are shown for different threshold
        values
        $\lambda$.}
    \label{fig:convergence_complex}
\end{figure}
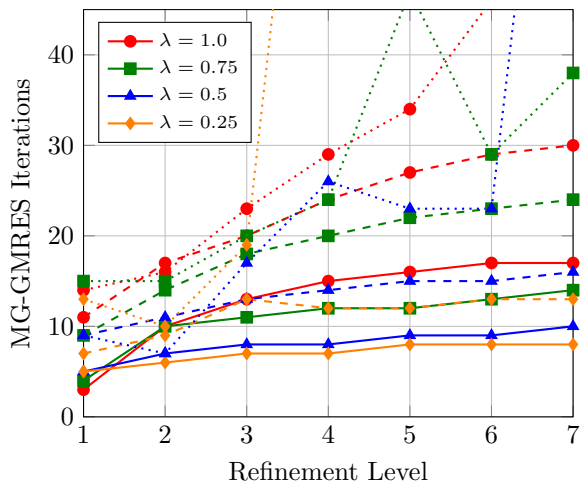

We investigate the performance of the multigrid preconditioner for the DG-SBM method on this complex deformed 2D
domain. To ensure a robust evaluation across different polynomial degrees for this complex geometry, we employ a full
$hp$-multigrid strategy. Figure~\ref{fig:convergence_complex} displays the number of GMRES iterations needed to reduce
the relative residual by a factor of $10^{-12}$. The figure illustrates the solver's behavior under mesh refinement for
various polynomial degrees and $\lambda$ values in 2D.

For polynomial degrees $p=1$ (solid lines) and $p=2$ (dashed lines), the multigrid preconditioner demonstrates good
scalability with mesh refinement. For $p=1$, iteration counts exhibit only a mild increase with mesh refinement across
all tested $\lambda$ values. For $\lambda=0.5$ and $\lambda=0.25$, the iteration counts remain particularly low and
stable, generally below 10 iterations even at the finest refinement level (level 7). For $p=2$, iteration counts are
higher but still show reasonable scalability, with $\lambda=0.5$ and $\lambda=0.25$ performing best, staying below 16
and 13 iterations respectively. For $p=3$ (dotted lines), the performance is more sensitive to $\lambda$. While
$\lambda=0.75$ (green dotted line) shows iteration counts increasing from 15 to 38, and $\lambda=0.5$ (blue dotted
line) increases from 9 to a failure at level 7 (100 iterations), other values show more significant growth or earlier
failure. Specifically, for $\lambda=1.0$ (red dotted line), iterations grow rapidly, reaching 100 (failure) at level 7.
For $\lambda=0.25$ (orange dotted line), the solver fails at level 4. These results for $p=3$ suggest that while the
preconditioner can be effective for some $\lambda$ values, its consistent performance for higher polynomial degrees on
complex geometries is more dependent on the choice of $\lambda$.

Overall, the results on the complex deformed 2D domain indicate that the DG-SBM multigrid preconditioner can achieve
efficient convergence, with iteration counts showing good scalability with mesh refinement for $p=1$ and $p=2$. The
choice of the threshold parameter $\lambda$ plays a significant role, especially for higher polynomial degrees,
influencing both the accuracy (as seen in Figure \ref{fig:solution_complex}) and the iterative solver performance.

In summary, the numerical investigations demonstrate that the proposed DG-SBM multigrid preconditioner generally
achieves good $h$-independent or mildly $h$-dependent convergence for linear ($p=1$) and quadratic ($p=2$) elements in
2D for both the unit sphere and the complex geometry. The threshold parameter $\lambda$ plays a crucial role. While
lower $\lambda$ values typically lead to lower $L_2$ errors and often improve solver stability for $p=1$, this trend
does not consistently hold for higher polynomial degrees. For $p \ge 2$, smaller $\lambda$ values, even if beneficial
for accuracy, can sometimes degrade solver performance, as observed in several test cases. Significant challenges
emerge for cubic ($p=3$) elements, particularly for complex geometries, where the current cell-wise SSOR smoother
exhibits reduced effectiveness, leading to increased iteration counts or solver failure.

\subsection{Comparison with DG-CutFEM}
In this section, we compare the performance of the proposed DG-SBM preconditioner with a Discontinuous Galerkin CutFEM
(DG-CutFEM) formulation. Both methods are evaluated using linear ($p=1$) and quadratic ($p=2$) elements. For the DG-SBM
computations, we utilize a threshold parameter of $\lambda=0.75$. To ensure a fair comparison, both solvers employ a
block-SSOR smoother with 3 smoothing steps. It is worth noting that Bergbauer~\cite{bergbauer2024high} also employed a
cell-wise smoother for DG-CutFEM discretizations. However, their approach\footnote{The specific order of the Chebyshev
    polynomial was confirmed through personal communication with the author.} utilized an additive cell-wise smoother
wrapped in a Chebyshev polynomial operator of order 10.

Table~\ref{tab:comparison_cutfem} presents the number of GMRES iterations required to reduce the residual by a factor
of $10^{-12}$ for a sequence of mesh refinements. The results indicate that the DG-SBM preconditioner exhibits
competitive performance compared to DG-CutFEM. For $p=1$, DG-SBM shows iteration counts very similar to DG-CutFEM, with
both methods demonstrating low mesh dependence. For $p=2$, DG-SBM also performs comparably, maintaining stable
iteration counts that are close to those observed for DG-CutFEM.

\begin{table}[h!]
    \centering
    \caption{Comparison of GMRES iteration counts for DG-SBM ($\lambda=0.75$), DG-CutFEM, and AMG applied to continuous SBM. DG-SBM and DG-CutFEM use a block-SSOR smoother with 3 smoothing steps.}
    \label{tab:comparison_cutfem}
    \begin{tabular}{ll|ccccccccc}
        \hline
        Method                      & Degree & 1   & 2   & 3   & 4   & 5   & 6   & 7   & 8   & 9   \\
        \hline
        \multirow{2}{*}{DG-SBM}     & $1$    & 4   & 6   & 9   & 10  & 11  & 13  & 13  & 14  & 15  \\
                                    & $2$    & 7   & 9   & 12  & 15  & 17  & 19  & 20  & 21  & 23  \\
        \cline{2-11}
        \multirow{2}{*}{DG-CutFEM}  & $1$    & 4   & 7   & 8   & 9   & 9   & 11  & 11  & 15  & 15  \\
                                    & $2$    & 8   & 11  & 14  & 17  & 16  & 19  & 19  & 24  & 24  \\
        \cline{2-11}
        \multirow{2}{*}{AMG + SBM } & $1$    & 2   & 2   & 1   & 11  & 12  & 13  & 14  & 16  & 16  \\
                                    & $2$    & --- & --- & --- & --- & --- & --- & --- & --- & --- \\
        \hline
    \end{tabular}
\end{table}

The results indicate that the DG-SBM preconditioner exhibits robustness comparable to the ghost-penalty stabilized
DG-CutFEM. Crucially, while DG-CutFEM requires ghost penalty stabilization and specialized quadrature on cut cells,
DG-SBM achieves similar iteration counts using only standard quadrature on background cells. This suggests that for
$p=1$ and $p=2$, SBM offers a more computationally efficient route to the same solver performance.

Furthermore, we include results for an algebraic multigrid (AMG) preconditioner applied to a standard continuous
Galerkin SBM formulation (AMG + SBM). For linear elements ($p=1$), AMG shows very low iteration counts on coarse
meshes, but the number of iterations increases and stabilizes as the mesh is refined, eventually exceeding those of the
DG-SBM and DG-CutFEM methods. More importantly, for quadratic elements ($p=2$), the AMG solver failed to converge, as
indicated by the dashes in Table~\ref{tab:comparison_cutfem}. This failure underscores the difficulty of using standard
algebraic multigrid methods for higher-order SBM discretizations, where the localized perturbations at the boundary can
significantly affect the spectral properties of the system matrix, making it challenging for AMG to construct effective
coarse-grid operators.

\section{Conclusion}
\label{sec:conclusion}
We have developed a geometric multigrid preconditioner for a Discontinuous Galerkin (DG) formulation of the Shifted
Boundary Method (SBM), aimed at efficiently solving the Poisson equation on domains with complex boundaries. Numerical
results demonstrate that the method, employing a cell-wise SSOR smoother within an $hp$-multigrid framework, achieves
good scalability and efficiency for linear ($p=1$) and quadratic ($p=2$) polynomial degrees in 2D.

The choice of the threshold parameter $\lambda$, which defines the set of active cells, significantly impacts solver
performance. This sensitivity to $\lambda$ is particularly evident for $p=2$ elements in 3D on complex geometries,
where solver stability can be critically affected. The influence of $\lambda$ on iteration counts is nuanced: while
smaller $\lambda$ values generally benefit solver performance for linear elements, this trend does not consistently
hold for higher polynomial degrees ($p \ge 2$). For these cases, the optimal $\lambda$ for solver efficiency may
differ, and smaller $\lambda$ values, even if they improve accuracy, can sometimes degrade solver performance.

Our initial 1D analysis indicated that SBM formulations can lead to system matrices with complex eigenvalues,
particularly when the true boundary lies within the surrogate cell (negative shifts). This insight is consistent with
numerical findings where smaller threshold values $\lambda$ (e.g., $\lambda=0.25$), which can increase the likelihood
of negative shifts, sometimes led to degraded convergence for higher-order elements.

The current cell-wise SSOR smoother shows limitations for higher polynomial degrees, particularly for $p \ge 3$ in 3D,
where iteration counts increase substantially or the solver fails to converge. This indicates that the smoothing
properties of the cell-wise SSOR are not sufficient for these more challenging SBM discretizations, possibly due to the
aforementioned spectral properties becoming more pronounced.

Future work should focus on addressing the limitations observed with the cell-wise SSOR smoother for higher polynomial
degrees ($p \ge 3$), particularly in 3D. The development and analysis of more advanced smoothers tailored to the
specific properties of DG-SBM systems are crucial. Investigating alternative smoothing strategies, such as Block Jacobi
(which was explored, though consistent significant improvements were not observed in the current study without further
tuning) potentially with Chebyshev acceleration, could be investigated. A thorough theoretical analysis of the
smoother's properties and its interaction with the SBM formulation is also warranted to guide the design of more
effective multigrid components. Furthermore, for higher polynomial degrees, especially when dealing with larger shift
distances, modifications to the SBM formulation itself might be necessary to enhance stability and improve smoother
performance.

\section*{Acknowledgments}
The author would like to thank Guglielmo Scovazzi and Guido Kanschat for insightful discussions, and Luca Heltai for
encouraging the exploration of this method.

\section*{Declarations}
The author declares support of two local feline agents (Micro and Conda) alongside language models (Gemini, Claude, ChatGPT) during text drafting. The final manuscript was audited by the author, who retains full accountability for all scientific content.

\FloatBarrier
\bibliographystyle{siam} 
\bibliography{literature} 
\end{document}

%% file: preamble.tex
\usepackage{amstext}
\usepackage{stmaryrd}
\usepackage{algpseudocode}

\usepackage[normalem]{ulem}
\usepackage{graphicx}

\usepackage{array}
\usepackage{verbatim}
\usepackage{booktabs}
\usepackage{calc}
\usepackage{textcomp}
\usepackage{multirow}

\usepackage[all]{xy}

\usepackage{cancel}
\usepackage{mathtools}
\usepackage{placeins}
\usepackage{xspace}

\usepackage{tikz,tikzscale}
\usetikzlibrary{calc}
\usetikzlibrary{shadings}

\usepackage[mode=multiuser,status=draft]{fixme} 
\fxsetup{envlayout=color}
\fxsetup{innerlayout=inline}
\fxsetup{targetlayout=colorcb}
\FXProvidesTargetLayout{color}

\fxusetheme{color}
\FXRegisterAuthor{ss}{envss}{SS}
\FXRegisterAuthor{mr}{envmr}{MR}
\FXRegisterAuthor{mw}{envmw}{MW}
\FXRegisterAuthor{mk}{envmk}{MK}

\usepackage{subcaption}
\usepackage[format=plain,indention=.5cm]{caption}
\usepackage{pgfplots}
\pgfplotsset{compat=1.18}
\usepgfplotslibrary{fillbetween}

\usepackage{bm,upgreek}

\newcommand{\DG}{{\mathbb{V}_h^{\text{DG}}}}
\newcommand{\DGell}{{\mathbb{V}_\ell^{\text{DG}}}}
\usepackage{graphicx} 
\usepackage[linesnumbered,vlined,commentsnumbered,ruled]{algorithm2e}

\definecolor{darkgreen}{rgb}{0.0, 0.5, 0.0}

